\documentclass[journal]{IEEEtran}

\usepackage{cite}
\usepackage{color}

\ifCLASSINFOpdf
   \usepackage[pdftex]{graphicx}
\else
\fi

\usepackage{amsmath}

\usepackage{array}

\begin{document}

\title{Computation of Convex Hull Prices \\ in Electricity Markets with Non-Convexities \\using Dantzig-Wolfe Decomposition}

\author{Panagiotis~Andrianesis,~\IEEEmembership{Member,~IEEE}, Dimitris Bertsimas, Michael C. Caramanis,~\IEEEmembership{Senior~Member,~IEEE}, and~William W. Hogan
\thanks{P. Andrianesis and M. C. Caramanis are with Boston University, Boston, MA: panosa@bu.edu, mcaraman@bu.edu. D. Bertsimas is with MIT, Cambridge, MA: dbertsim@mit.edu. W. W. Hogan is with Harvard University, Cambridge, MA, william\_hogan@hks.harvard.edu.
Research partially supported by NSF AitF 1733827.
}
}

\maketitle
\begin{abstract}
The presence of non-convexities in electricity markets has been an active research area for about two decades.
The --- inevitable under current marginal cost pricing --- problem 
of guaranteeing that no market participant incurs losses in the day-ahead market 
is addressed in current practice through make-whole payments \emph{a.k.a.} uplift. 
Alternative pricing rules have been studied to deal with this problem. 
Among them, Convex Hull (CH) prices associated with minimum uplift have attracted significant attention.
Several US Independent System Operators (ISOs) have considered CH prices but resorted to approximations, 
mainly because determining exact CH prices is computationally challenging, while providing little intuition about the price formation rational{\color{black}e}.
In this paper, we describe the CH price estimation problem 
by relying on Dantzig-Wolfe decomposition and Column Generation, 
{\color{black}as a tractable, highly paralellizable, and exact method 
--- i.e., yielding exact, not approximate, CH prices --- 
with guaranteed finite convergence.}
Moreover, the approach provides intuition on the underlying price formation rational{\color{black}e}. 
A test bed of stylized examples provide an exposition of the intuition in the CH price formation. 
In addition, a realistic ISO dataset is used to support scalability and validate the proof-of-concept.
\end{abstract}


\IEEEpeerreviewmaketitle

\section{Introduction}

\IEEEPARstart{M}{arginal} cost pricing based on spot pricing under convexity assumptions  \cite{CaramanisEtAl_1982, SchweppeEtAl_1988} has been the standard practice in organized electricity markets.
However, in the presence of non-convexities 
(mainly due to unit commitment costs and technical constraints, e.g., minimum output requirements),
marginal cost pricing from a restricted convex subproblem cannot guarantee support of the solution where 
market participants recover their as-bid production costs.  There may be no market-clearing prices for the economically efficient solution.
This problem has been typically dealt with the provision of ``uplift'' side-payments to market participants, 
to make them whole.
Over the last two decades, motivated mainly by the electricity market paradigm, 
pricing in markets with non-convexities has attracted significant attention.
Approaches have ranged from standard marginal cost pricing with recovery mechanisms --- e.g., \cite{ONeillEtAl_2005, AndrianesisEtAl_2013a, AndrianesisEtAl_2013b}, 
to mechanism designs that ``minimize,'' internalize or, sometimes, even eliminate uplifts --- e.g., \cite{HoganRing_2003, GribikEtAl_2007, LibAnd_2016, ONeillEtAl_2017}.
A critical review of pricing rules employed in markets with non-convex costs is provided in \cite{LibAnd_2016}.

Recently, FERC (Federal Energy Regulatory Commission) 
initiated a discussion on price formation \cite{FERC_2014},
with several US Independent System Operators (ISOs) exploring Convex Hull (CH) prices,  
--- first suggested in \cite{GribikEtAl_2007} ---
and their approximations, 
often called Extended Locational Marginal Prices (ELMPs) 
\cite{FERC_2014, WangEtAl_2016, ChenWang_2017, PJM_2018, MISO_2019}.
CH pricing was introduced in \cite{GribikEtAl_2007} as the CH of the aggregate cost function, i.e., 
the convex function that is closest to approximating the aggregate cost function from below.
It is equivalently derived by Lagrangian dualization or CH relaxation, i.e.,
considering the CH of individual components, as pointed out already in \cite{GribikEtAl_2007} 
--- see also the interesting discussion in \cite{Chao_2019}. 

Interestingly, up until the early 2000s, 
the Unit Commitment (UC) problem itself was traditionally solved using Lagrangian Relaxation (LR),
first by vertically integrated utilities, 
and later by ISOs that ``inherited''
LR solvers --- possibly the only commercial option at the time --- in the first market implementations.
Taking advantage of Mixed Integer Linear Programming (MILP) solver advances, 
in 2005,
PJM replaced LR with MILP achieving cost savings of about \$500M/year; 
by 2017, all US ISOs switched to MILP with estimated cost savings of more than \$2B/year \cite{ONeill_2020}.
Currently, the Day-Ahead (DA) market is based on a sequence of Security Constrained UC (SCUC) and Security Constrained Economic Dispatch (SCED), followed after closure of the DA Market by an ISO executed Reliability Unit Commitment (RUC) process.
{\color{black}Despite the advancements in the literature of stochastic and particularly robust optimization \cite{bertsim} 
in both DA and real-time settings (e.g., \cite{ding1, ding2}), 
the UC problem is still being solved in practice as a deterministic problem.}
DA prices are finally obtained by a pricing run with fixed unit commitments, 
either before (e.g., PJM) or after (e.g., NYISO) the RUC process.
Some ISOs proceeded to ELMPs, 
by relaxing integrality constraints for fast-start units and using ELMPs as a proxy to CH prices \cite{WangEtAl_2016}.
Integer Relaxation (IR) of the MILP UC problem is considered in \cite{Chao_2019} as a simple method that
produces good approximations 
--- sometimes even exact --- to the (equivalent) Lagrangian Dual (LD) and CH relaxations,
which, although referred to as ``an ideal solutions,'' can remain computationally prohibitive.

Several works have attempted to overcome the computational difficulties by relying on two main approaches.
An early approach applied sub-gradient methods \cite{WangEtAl_2009, WangEtAl_2013, WangEtAl_2013a, WangEtAl_2013b, ItoEtAl_2013}, 
whereas a later and methods focused on identifying the CH of individual generators 
through convex primal formulations \cite{HuaBaldick_2017} 
---  also including an AC Optimal Power Flow setting \cite{GarciaEtAl_2020}, 
extended formulations \cite{YangEtAl_2019, ChenEtAl_2020, YuEtAl_2020},
a network reformulation \cite{AlvarezEtAl_2020},  
and Benders decomposition leveraging advances in thermal generator CH formulations \cite{KnuevenEtAl_2019}.
Despite the aforementioned efforts,
two important barriers to CH price implementation remain: 
(i) computational challenges, 
and (ii) opacity of their properties \cite{PJM_2018}.
An insightful work on the latter, \cite{SchiroEtAl_2016}, 
uses representative stylized examples to illustrate some arguably counter intuitive CH price properties.

Our aim here is to address the computational challenge while providing intuition for the CH price formation.
Our approach is inspired by related experience in solving large-scale optimization problems in other application domains, in particular crew scheduling \cite{AndKoz_2014}. 
The Operations Research theory underpinnings date back to the 60's 
and the seminal work of Dantzig and Wolfe \cite{DantzigWolfe_1960} on Generalized Linear Programming (LP), 
\emph{a.k.a.} Dantzig-Wolfe (D-W) decomposition or Column Generation (CG).
D-W decomposition can be also viewed as a problem characterization 
rendering the problem suitable to a CG algorithm.
The relationship of Generalized LP 
as a solution method for the LD dates also back to the 70's 
--- see e.g., \cite{Geoffrion_1974, MagnantiEtAl_1976}.
Around the same time, a CG process was also proposed for approximating competitive equilibria in a piecewise linear economy \cite{ManneEtAl_1980}.
Later, around the 90's, the work of \cite{BarnhartEtAl_1998} on ``branch-and-price'' --- among others \cite{Vanderbeck_2000} --- 
positioned CG as a powerful solution method for huge integer problems, 
with two common illustrative applications: the generalized assignment problem and crew scheduling.
Since then, CG has been successfully applied to large scale integer programming \cite{LubbeckeDesrosiers_2005}, 
and is possibly the only commercially available option for scheduling problems in the airline industry.

Our analysis draws from CH price ``first principles'' \cite{GribikEtAl_2007}, and our main focus has two elements.
First, 
we provide a D-W characterization of the UC problem, without changing any of the original functions or constraints, 
whose LP relaxation is equivalent to the solution of the LD of the original MILP formulation;
we then present a CG algorithm to solve the LP relaxation, 
and derive the CH prices.
Second, we illustrate the applicability of the approach to 
(i) stylized examples in \cite{SchiroEtAl_2016} that provide intuition into the CH price formation,
(ii) a more detailed ramp-constrained example from \cite{ChenEtAl_2020}, 
to show that our approach can handle features that could not be practically addressed without reformulation of the UC problem,
and (iii) a large-scale ISO dataset \cite{FERC_RTO} with about 1000 generators, 
thus illustrating the potential scalability.

{\color{black} Application of D-W characterization and a CG solution algorithm to the CH pricing problem provides a new direction, 
which seems as a natural fit,
with several comparative advantages. 
It proposes a tractable, highly parallelizable, and exact method 
--- i.e., a method that yields exact, not approximate, CH prices --- 
with guaranteed finite convergence.
It also presents a unique economic interpretation and intuition about the CH price formation.
This is not a coincidence; 
it originates from the economic interpretation of the D-W decomposition 
and its link with the Lagrangian Dual.
This work identifies the relationship of D-W decomposition to the derivation of CH prices, 
while aligning with the essence of CH prices 
and their property of supporting the market solution with ``minimum uplift.''
Notably, the proposed approach moves away from complicated extended formulations, 
and can cope directly with existing MILP implementations of unit characteristics and constraints.
It shapes the CH, with a simple ``convexity constraint,'' 
while moving the complexity of the MILP-constrained units to small sub-problems.
The derivation of the convex hull directly from the convexity constraint of the D-W characterization, 
without resorting to special case unit specific formulations, 
is a unique contribution to the so far known literature.
As such, the proposed approach is generalizable and can accommodate any change in the unit formulations, 
without requiring approximations/simplifications.
Last but not least, the preliminary computational experiments in a large-scale dataset are particularly encouraging, with a high potential for practical implementation.}

The remainder of the paper is organized as follows.
In Section \ref{CHdescr}, we define CH prices for a stylized UC problem.
In Section \ref{DWref}, we present a D-W reformulation of the problem and sketch the CG solution algorithm.
In Section \ref{Stylized}, we detail the application to stylized examples.
In Section \ref{MoreTests}, we discuss {\color{black}the comparative advantages of the proposed method and} demonstrate the computational tractability on realistic test cases.
Lastly, Section \ref{Conclusions} concludes and provides directions for further research.

\section{Description of CH Prices} \label{CHdescr}

To fix ideas, 
but without loss of generality, 
we assume that $x_{i,t}$ refers to power output (energy for a dispatch interval) 
and $y_{i,t}$ refers to the discrete variables (e.g., the on/off status) of unit $i$, at time period $t$.
Let $\mathcal{I}$ denote the set of generation units,
and $\mathcal{T} = \{1, 2, ..., T\}$ the set of time periods, 
where $T$ is the length of the optimization horizon.
In what follows, we refer to $t$ as an hour for simplicity.
For brevity, and occasionally with some abuse of notation, 
we use $\mathbf{x}_i, \mathbf{y}_i$,   
to denote vectors of unit $i$ comprising the respective variables, $x_{i,t}, y_{i,t}$, $ \forall t \in \mathcal{T}$;
at a more abstract level, we use $\mathbf{x}$ and  $\mathbf{y}$ as vectors comprising the respective variables, $x_{i,t}$ and $y_{i,t}$, 
$\forall i \in \mathcal{I}$, and $ \forall t \in \mathcal{T}$.
Let also $D_t$ denote the demand for energy at hour $t$. 
A basic UC problem is formulated as follows:
\begin{subequations} \label{UCprob}
\begin{align} 
\underset{\mathbf{x}, \mathbf{y}}{\min} \,  f(\mathbf{x},\mathbf{y}) & = \sum_{i \in \mathcal{I}} f_i(\mathbf{x}_i, \mathbf{y}_i), &    \label{UCobj} \\
\text{subject to: }\qquad  \sum_{i \in \mathcal{I}} x_{i,t} & = D_t,  & \forall t \in \mathcal{T},  \label{UCpbal} \\
(\mathbf{x}_{i}, \mathbf{y}_{i}) & \in \mathcal{Z}_i,  & \forall i \in \mathcal{I}.  \label{UCunitcon}
\end{align}
\end{subequations}

The objective function \eqref{UCobj} minimizes the aggregate commitment and dispatch costs for all units, 
with $f_i (\mathbf{x}_{i}, \mathbf{y}_{i}) $ denoting the respective cost for unit $i$ 
over the entire optimization horizon,
{\color{black}typically using a mixed integer linear representation.}
Constraints \eqref{UCpbal} represent the power balance for all hours;
they can be straightforwardly extended to include all system constraints 
(e.g., power balance, reserve requirements, etc., 
and these can be extended to incorporate linear representations of transmission constraints.)
{\color{black}Constraints \eqref{UCunitcon} represent all applicable unit specific constraints 
(capacity, ramps, etc.) that are present in actual market implementations,
typically using a mixed integer linear representation,}
as well as integrality constraints for $\mathbf{y}$ variables,
with $\mathcal{Z}_i$ denoting the set of constraints for unit $i$ for the entire horizon.
We use $\mathcal{Z}_i$ to ease the exposition for reasons that will soon become apparent.
In this respect, constraints \eqref{UCunitcon} define a \emph{feasible schedule} for unit $i$.
It may not be evident at this point, but this last remark is arguably the most important and key to our analysis.
We also note that, although throughout the paper we use for simplicity the term UC, 
we essentially refer to a UC and Economic Dispatch problem in the broad sense.
Our analysis can also cater to current ISO needs and practices (e.g., SCUC, SCED, RUC).

Standard marginal cost pricing relates to the dual variable of constraint \eqref{UCpbal},
after fixing discrete variables $\mathbf{y}$ to their optimal values.
Enforcing such explicit equality constraints 
is used in \cite{ONeillEtAl_2005} to price integralities, 
under standard marginal cost pricing, 
and immediately derive the make-whole payments, 
in what has been called the ``IP pricing,'' --- see also \cite{LibAnd_2016}.
ELMPs are currently derived by relaxing integer variables of only the fast-start units.
In general, the IR of the MILP UC problem \eqref{UCprob} yields less tight convex approximations.

\begin{subequations} \label{LagDual}
Revisiting \cite{GribikEtAl_2007}, 
there are two equivalent methods to derive CH prices.
The first method prices out complicating system constraints through Lagrangian dualization. 
It is based on dualizing constraint \eqref{UCpbal}, 
with $\lambda_t$ the respective dual and $\boldsymbol \lambda$ an appropriate vector, 
and with CH prices obtained by the solution of the following problem, which we denote by LD:
\begin{align} 
\underset{\boldsymbol{\lambda}}{\max} \,
q(\boldsymbol{\lambda}), \label{LagDual1} 
\end{align}
where
$q(\boldsymbol{\lambda})$ is the dual function of the Lagrangian, $\mathcal{L}(\mathbf{x}, \mathbf{y}, \boldsymbol{\lambda})$,
\begin{equation} \label{qDual}
 q(\boldsymbol{\lambda}) =
\underset{(\mathbf{x}_{i}, \mathbf{y}_{i}) \in \mathcal{Z}_i, \, \forall i \in \mathcal{I}}{\inf} \mathcal{L}(\mathbf{x}, \mathbf{y}, \boldsymbol{\lambda}),  
\end{equation}
\begin{equation} \label{partLagr}
\mathcal{L}(\mathbf{x}, \mathbf{y}, \boldsymbol{\lambda})  =  
\sum_{i \in \mathcal{I}} f_i(\mathbf{x}_i, \mathbf{y}_i) - \sum_{t \in \mathcal{T}} \lambda_t \Big(\sum_{i\in \mathcal{I}} x_{i,t} - D_t \Big).
\end{equation}
\end{subequations}
The second method is based on describing the CH of individual components 
usually employing tight/extended formulations,
and the convex envelope of the objective function.
However, regardless of the derivation method, 
computational intractability persists, in particular for general UC settings and unit characteristics.

{\color{black}CH prices are named after the CH of the feasible sets ($\mathcal{Z}_i$), 
and should be distinguished from the CH 
of the entire UC MILP problem \eqref{UCprob}.
Equivalently, CH prices could be named after the LD problem, 
and referred to as LD prices.
It is important to note that  
CH (or equivalently LD) prices solve the LD problem,
i.e., a pricing problem, 
which is not the same as the UC MILP problem (there is a duality gap).
In the pricing problem 
(which is solved during the pricing run in actual markets, 
as we also mentioned in our Introduction),
we are only interested in obtaining the prices 
(i.e., the $\boldsymbol{\lambda}$'s)
and not the quantities 
(i.e., the $\mathbf{x}$'s and $\mathbf{y}$'s).
The latter (i.e., the quantities) 
are obtained by the solution of the UC MILP problem.
CH (or LD) prices refer to the solution of the LD problem \eqref{LagDual}.
They are exact if they exactly solve the LD problem,
which also means that they exactly identify the CH of the feasible sets.
They are approximate if they approximately solve the LD problem, 
which also means that they approximately identify the CH of the feasible sets.
Current implementations of ELMPs are not exact CH prices.
They are often considered as CH approximations,
but in reality they are limited implementations of IR.}

CH prices are also often called ``minimum uplift'' prices.
This characterization
essentially refers to the minimization of the total Lost Opportunity Costs (LOCs) 
in addition to what is referred to as ``Product Revenue Shortfall'' (PRS) in \cite{SchiroEtAl_2016}.
For brevity, we focus our discussion on LOCs, and refer the reader to \cite{SchiroEtAl_2016} for PRS.  
Notably, ``uplift minimization'' is directly related to the duality gap 
between {\color{black}the UC MILP problem} \eqref{UCprob} and {\color{black}the LD problem} \eqref{LagDual}.

Compensating for LOCs, even for units that are offline as a result of the UC solution, 
is a means to deter inefficient self-scheduling.
More specifically, given prices $\boldsymbol{\lambda}$, 
let us denote the profit of unit $i$ by $\phi_i(\mathbf{x}_i, \mathbf{y}_i; \boldsymbol{\lambda})$, where:
\begin{equation} \label{profit}
 \phi_i(\mathbf{x}_i, \mathbf{y}_i;\boldsymbol{\lambda}) = \sum_{t \in \mathcal{T}} \lambda_t  x_{i,t} -  f_i(\mathbf{x}_i, \mathbf{y}_i).   
\end{equation}
Unit $i$, following the solution of \eqref{UCprob},
assuming a market schedule $(\mathbf{x}_i^M, \mathbf{y}_i^M)$ and prices $\boldsymbol{\lambda}$,
would gain $\phi_i (\mathbf{x}_i^M, \mathbf{y}_i^M; \boldsymbol{\lambda})$.
However, given prices $\boldsymbol{\lambda}$, 
had the unit optimally self-scheduled, 
i.e., had the unit solved the following profit maximization problem: 
\begin{equation} \label{Self-Scheduling}
(\mathbf{x}_{i}^S, \mathbf{y}_{i}^S; \boldsymbol{\lambda}) \in \underset{ (\mathbf{x}_{i}, \mathbf{y}_{i}) \in \mathcal{Z}_i}{\text{argmax}} \,\,  \phi_i (\mathbf{x}_i, \mathbf{y}_i; \boldsymbol{\lambda}),
\end{equation}
the tentative profit would have been $\phi_i (\mathbf{x}_i^S, \mathbf{y}_i^S; \boldsymbol{\lambda})$.
Hence, the LOC for unit $i$, LOC$_i$, is given by:
\begin{equation} \label{LOC}
\text{LOC}_i = \phi_i (\mathbf{x}_i^S, \mathbf{y}_i^S; \boldsymbol{\lambda}) - \phi_i (\mathbf{x}_i^M, \mathbf{y}_i^M; \boldsymbol{\lambda}).
\end{equation}
CH prices are complemented by minimal LOCs according to \eqref{LOC} plus the PRS, to support the market solution by making the unit indifferent between following the market schedule $(\mathbf{x}_i^M, \mathbf{y}_i^M)$ and self-scheduling $(\mathbf{x}_i^S, \mathbf{y}_i^S)$.

\section{D-W Characterization and CG Algorithm} \label{DWref}

{\color{black}
Before proceeding with the D-W characterization, 
let us first note, in a parenthesis, that stylized formulations of the UC problem would describe it as follows:
\begin{equation*}
\underset{\mathbf{x}}{\min} \, \tilde{f}(\mathbf{x}) = \sum_{i \in \mathcal{I}} \tilde{f}_i ( \mathbf{x}_i ), \, \text{subject to: } \eqref{UCpbal},  \mathbf{x}_i \in \mathcal{X}_i,  \, \forall i \in \mathcal{I},
\end{equation*}
considering $\mathbf{x}_i$ as a feasible dispatch schedule, 
internalizing the discrete (status) variables ($\mathbf{y}_i$) in the feasibility constraint set $\mathcal{X}_i$,
and defining the cost function $\tilde{f}_i(\cdot)$ so as to account for both dispatch and commitment costs.}
Indeed, for a certain dispatch schedule, the unit status and cost can be considered as endogenously determined --- see for instance \cite{FeltKiwiel_2001, Lemarechal_2007, LunaEtAl_2020} that relate to an older work \cite{BertsekasEtAl_1983} for such formulations --- among others.
The main reason for this parenthesis is that the above formulation --- which is equivalent to \eqref{UCprob},
provides a natural interpretation of the proposed D-W characterization of the CH.

{\color{black}
Let us return to the UC problem formulation as presented in \eqref{UCprob}.
The state-of-the-art commercial products, 
i.e., the Market Management Systems that are in use in actual markets, employ a MILP formulation.}
{\color{black}In the state-of-the-art implementations,}  
{\color{black}the objective function $f(\cdot)$ is a mixed integer linear representation of the aggregate commitment and dispatch costs described by $f_i(\cdot)$ 
for each unit $i$.
Similarly, the set of constraints $\mathcal{Z}_i$ for unit $i$ is represented by mixed integer linear constraints.
The problem is solved using commercially available MILP solvers,
which allow computational tractability.
For obtaining CH prices, 
we make no changes and no assumptions to the form of the constraints that characterize the units.
We do not need to apply any special treatment for specific type of constraints (e.g., ramping).
We use exactly the same formulations that are currently implemented in the state-of-the-art UC implementations.
Our analysis is based on these facts,
{\color{black}(i.e., assumes the commercial state-of-the-art mixed integer linear representations of objective function and constraints)}
and requires that the set $\mathcal{Z}_i$ is bounded, 
which holds for all units.}

Consider an equivalent, {\color{black} albeit arguably ``unusual,''} formulation of the UC problem \eqref{UCprob},
using variable $z_i$ to describe a feasible schedule of unit $i$, 
for the entire optimization horizon, i.e., 
$z_i := (\mathbf{x}_{i}, \mathbf{y}_{i}) \in \mathcal{Z}_i$.
{\color{black} 
Let $n_i \in \mathcal{N}_i$ be an index of unit $i$ feasible schedules, 
    which are contained in the set $\mathcal{N}_i$.}\footnote{
{\color{black} We clarify that we do not need an assumption for the number of schedules to be finite.
It suffices that the feasible sets $\mathcal{Z}_i$ are bounded \cite{MagnantiEtAl_1976}.
}
}
Considering the UC formulation \eqref{UCprob},
each variable $z_{i}^{[n_i]}$ corresponds to a specific feasible schedule, 
denoted by $( \mathbf{\hat x}_{i}^{[n_i]}, \mathbf{\hat y}_{i}^{[n_i]})$ 
and has a cost denoted by $\hat c_i^{[n_i]} = f_i( \mathbf{ \hat x}_{i}^{[n_i]}, \mathbf{\hat y}_{i}^{[n_i]})$.
Equivalently, a feasible schedule is also determined by $\mathbf{\hat x}_{i}^{[n_i]}$ and $\hat c_i^{[n_i]}$.
We will use the two representations interchangeably to facilitate the exposition.

\begin{subequations} \label{UC-DW}
Following the above, the UC problem can be expressed as:
\begin{align} 
\underset{\mathbf{z}}{\min} \,  \sum_{ i \in \mathcal{I}, n_i \in \mathcal{N}_i } {\color{black} f_i( \mathbf{ \hat x}_{i}^{[n_i]}, \mathbf{\hat y}_{i}^{[n_i]}) } \, z_{i}^{[n_i]}, &  \label{UCobj_z} \\
\text{subject to:} \, \sum_{i \in \mathcal{I}, n_i \in \mathcal{N}_i} \hat x_{i,t}^{[n_i]} z_{i}^{[n_i]}  = D_t, \,\,& \forall t \in \mathcal{T},  \label{UCpbal_z}\\ 
 \sum_{n_i \in \mathcal{N}_i} z_{i}^{[n_i]}  = 1, \quad & \forall i \in \mathcal{I}, \label{UCunit_z1}\\
 z_{i}^{[n_i]} \in \{0, 1\}, \quad  & \forall i  \in \mathcal{I}, n_i \in \mathcal{N}_i.  \label{UCunit_z2}
\end{align}
\end{subequations} 
Solving \eqref{UC-DW}, 
accounting for all feasible schedules, 
is equivalent to solving \eqref{UCprob}. 
Objective function \eqref{UCobj_z} is equivalent to \eqref{UCobj}. 
Constraint \eqref{UCpbal_z} represents the power balance and is equivalent to \eqref{UCpbal}. 
{\color{black}Of note that $ f_i( \mathbf{ \hat x}_{i}^{[n_i]}, \mathbf{\hat y}_{i}^{[n_i]})$ and $ \hat x_{i,t}^{[n_i]}$ are parameters, and hence, both the objective function and the constraints are linear in $z_{i}^{[n_i]}$.}
Constraint \eqref{UCunit_z1} requires to select exactly one feasible schedule represented by the binary variable $z_{i}^{[n_i]}$.
{\color{black}Since the optimal solution of the UC problem \eqref{UCprob} --- say $(\mathbf{ \hat x}_{i}^{*}, \mathbf{\hat y}_{i}^{*})  \in \mathcal{Z}_i$ is itself one of the candidate feasible solutions --- denoted by $z_{i}^{*}$, then $z_{i}^{*}$ optimally solves problem \eqref{UC-DW}, and the optimal objective function value of \eqref{UCobj_z} will be equal to the respective value of \eqref{UCobj}.} 
Equivalently, problem \eqref{UC-DW} can be viewed as a set partitioning problem.

At a first glance, 
the effort of brute-force solving \eqref{UC-DW} would seem hopeless, 
as it would require not only an identification of all feasible schedules,
but also the solution of a huge integer problem  
--- \eqref{UC-DW} is an integer linear programming problem.
That said, we naturally associate \eqref{UC-DW} with D-W decomposition and CG --- see e.g., \cite{BarnhartEtAl_1998}.
However, solving a large UC problem with CG does not yet seem promising --- see e.g., a recent work \cite{Dupin_2019}.
Fortunately, we remind the reader that our goal is not to solve this integer problem,
but the LD \eqref{LagDual}.
With this in mind, we recall a result which dates back to the 70's \cite{Geoffrion_1974},
that the LP relaxation of \eqref{UC-DW}, 
solves the LD \eqref{LagDual}; 
this is also clearly stated in \cite{MagnantiEtAl_1976},
that Generalized LP --- \emph{a.k.a.} D-W decomposition --- solves the dual.
{\color{black}In other words, 
the LP relaxation of problem \eqref{UC-DW} is the D-W characterization of the UC problem 
--- equivalently stated as \eqref{UCprob} or \eqref{UC-DW} ---
and solves the LD of the UC problem, 
which prices out (dualizes) the system constraints.}
Indeed, we know that D-W decomposition or CG is essentially a cutting plane method applied to the LD \cite{Ber},
which we find as particularly instructive in our problem setting, 
hopefully enhancing the intuition into the formation of CH prices.

{\color{black}
For clarity, let us provide the LP relaxation of problem \eqref{UC-DW} below.
\begin{subequations} \label{UC-DW-LP}
\begin{align} 
 \underset{\mathbf{z}}{\min} \, g(\mathbf{z}) = \sum_{ i \in \mathcal{I}, n_i \in \mathcal{N}_i } f_i( \mathbf{ \hat x}_{i}^{[n_i]}, & \mathbf{\hat y}_{i}^{[n_i]}) \, z_{i}^{[n_i]}, \quad  \label{UCobj_z-LP} \\
\text{subject to:} \, \sum_{i \in \mathcal{I}, n_i \in \mathcal{N}_i} \hat x_{i,t}^{[n_i]} z_{i}^{[n_i]}  = D_t, \,\,& \forall t \in \mathcal{T},  \label{UCpbal_z-LP}\\ 
 \sum_{n_i \in \mathcal{N}_i} z_{i}^{[n_i]}  = 1, \quad & \forall i \in \mathcal{I}, \label{UCunit_z1-LP}\\
 z_{i}^{[n_i]} \geq 0, \quad  & \forall i  \in \mathcal{I}, n_i \in \mathcal{N}_i.  \label{UCunit_z2-LP}
\end{align}
\end{subequations} 
Note that limiting the upper bound to $z_{i}^{[n_i]} \leq 1$ is enforced by \eqref{UCunit_z1-LP}. 
Following \cite{MagnantiEtAl_1976}, problem \eqref{UC-DW-LP} is equivalent to the LD problem \eqref{LagDual}; the solution of \eqref{UC-DW-LP} also solves \eqref{LagDual}.
}

Many approaches to CH pricing involve restrictions on the form of the objective function \cite{Chao_2019} 
or reformulations of the UC problem such as dynamic programming characterizations \cite{KnuevenEtAl_2018, BacciEtAl_2019}.
The D-W approach starts with a natural formulation, which does not require any change in the UC functions or constraints.  
This allows a wide degree of flexibility in defining the UC problem.
This representation allows us to characterize the full CH \cite{Falk_1969}.
Essentially, the same formulation is applied in crew scheduling problems \cite{BarnhartEtAl_1998}. 
If we had to deal with the full listing of all the feasible points, 
the computational problem would be overwhelming.  
But the D-W method provides a natural CG technique 
that uses sub-problems to produce what is in practice a relatively small number of feasible points 
adequate to characterize the solution of the CH relaxation.
The intuition is that the same approach would work for the UC problem, 
without requiring any reformulations of the UC model, 
but with the computational feasibility observed in the related very large crew scheduling problems.

The solution of {\color{black} the LP problem \eqref{UC-DW-LP}, which is the LP relaxation of problem \eqref{UC-DW},}
employing a CG algorithm, proceeds as follows.

First, consider a subset of feasible schedules for each unit $i$, 
say $\mathcal{N}_i^{(k)} \subset \mathcal{N}_i$,
where $k$ is an iteration counter, i.e., 
$\mathcal{N}_i^{(k)}$ contains all feasible schedules of unit $i$ available at iteration $k$.
\begin{subequations} \label{UC_RMP}
Using these schedules, a Restricted Master Problem (RMP) at iteration $k$, RMP$^{(k)}$, is formulated as follows:
\begin{align}
\text{RMP}^{(k)}: \qquad \quad  \underset{\mathbf{z}}{\min} \, g^{(k)}(\mathbf{z})  = & \sum_{i \in \mathcal{I}, \, n_i \in \mathcal{N}_i^{(k)}} \hat c_i^{[n_i]} z_{i}^{[n_i]}, \label{UCobjRMP}   \\ 
\text{subject to:} \sum_{i \in \mathcal{I}, n_i \in \mathcal{N}_i^{(k)}} \hat x_{i,t}^{[n_i]} z_i^{[n_i]} & =  D_t, \forall t \in \mathcal{T}, \rightarrow \lambda_t^{(k)},  \label{UCpbalRMP} \\
    \sum_{n_i \in \mathcal{N}_i^{(k)}} z_{i}^{[n_i]} &= 1, \,\, \forall i \in \mathcal{I},  \rightarrow \pi_i^{(k)},  \label{UCunit_zRMP}
\end{align}
with $z_{i}^{[n_i]} \geq 0, \, \forall i \in \mathcal{I}, \, n_i \in \mathcal{N}_i^{(k)}$,
and $\lambda_t^{(k)}$, $\pi_i^{(k)}$, the duals of constraints \eqref{UCpbalRMP}, \eqref{UCunit_zRMP}, respectively.
\end{subequations} 
{\color{black}We clarify that problem \eqref{UC_RMP} is called ``Restricted'' because it contains only a subset of schedules compared with the ``Master Problem,'' as is often called problem \eqref{UC-DW-LP}}.

{\color{black}The formulation of the initial RMP in \eqref{UC_RMP} assumes feasibility.  
Indeed, 
this can be easily achieved by populating the initial RMP 
with a UC feasible solution,  
e.g., with the feasible schedules derived by the solution of \eqref{UCprob}, 
if such a solution is available at the time.
Alternatively, we can populate the initial RMP 
with a trivial feasible schedule for each unit, 
e.g., considering that the unit is offline.
To ensure that a feasible solution exists for problem \eqref{UC_RMP},
we can introduce slack variables in the constraints, with appropriate penalties,
in the same way we use slack variables in UC MILP formulations.
This is also standard practice in CG implementations.}\footnote{{\color{black}
We refer to the examples in the sections that follow 
for an illustration of the initialization and the use of slack variables.}} 

Then, CG ``generates'' new columns (feasible schedules) and adds them to the RMP, 
as long as they have a negative reduced cost.
The reduced cost $rc(\cdot)$ of feasible schedule $( \mathbf {x}_{i},  \mathbf{y}_{i} )$ is given by:
\begin{equation} \label{rc} 
 rc_i^{(k)}( \mathbf{x}_{i}, \mathbf{y}_{i})  =  f_i ( \mathbf{x}_{i}, \mathbf{y}_{i} ) - \sum_{t \in \mathcal{T}} \lambda_t^{(k)} x_{i,t} - \pi_i^{(k)}.
\end{equation}
Hence, at iteration $k$,
feasible schedules with potential negative reduced cost 
can be obtained by the solution of the following sub-problem, for each unit $i$:
\begin{subequations}\label{SubProblem}
\begin{align} 
\text{Sub}_i^{(k)}: \underset{\mathbf{x}_i, \mathbf{y}_i}{\min} \,  h_i^{(k)}(\mathbf{x}_i, \mathbf{y}_i) = f_i(\mathbf{x}_i, \mathbf{y}_i) - \sum_{t \in \mathcal{T}} \lambda_t^{(k)}  x_{i,t},\\
\text{subject to: } \qquad   (\mathbf{x}_{i}, \mathbf{y}_{i}) \in \mathcal{Z}_i, \qquad \qquad \qquad 
\end{align}
\end{subequations}
where we dropped the constant term $-\pi_i^{(k)}$ from the objective function.
If the schedule obtained by the solution of \eqref{SubProblem} has a negative reduced cost,
which is calculated using \eqref{rc} 
--- essentially by adding back the term $-\pi_i^{(k)} $ to the objective function value $h_i^{(k)}$,
then a new column is added to the RMP corresponding to this schedule.
The RMP is solved again, 
and the algorithm terminates when no new feasible schedule with negative reduced cost can be found.
{\color{black}Notably, 
a schedule that already exists in the RMP cannot have a negative reduced cost, 
hence, no feasible schedule can be repeatedly added in the RMP.}
Obviously, we will not always be able to find negative reduced cost schedules at each iteration for all units;
this has been the reason that we did not use superscript $(k)$ instead of $[n_i]$ for the feasible schedule.
{\color{black} Finite convergence of the CG algorithm is guaranteed 
due to the mixed integer linear representation of objective function $f_i(\cdot)$ and constraint sets $\mathcal{Z}_i$.}
{\color{black} Specifically, MILP-constrained feasible sets $\mathcal{Z}_i$ have a polyhedral convex hull, 
with finitely many extreme points.
The RMP can describe any feasible point 
as a convex combination of finite set of extreme points.
Hence, since each iteration --- solving MILP subproblems \eqref{SubProblem} --- includes at least one extreme point in the RMP, 
the method converges finitely.\footnote{
{\color{black}We refer the interested reader to \cite{MagnantiEtAl_1976} for convergence proofs of the method in more general settings than the one required for the solution of the CH pricing problem.}
}}

We further elaborate on the intuition that the CG algorithm provides.
Evidently, the RMP selects fractional feasible schedules, 
in fact a convex combination of them 
--- see constraint \eqref{UCunit_zRMP} that is often called a ``convexity constraint'' for a reason.
CG terminates by ``shaping'' the CH yielding the CH prices.
Indeed, had we asked a similar question as in \cite{BergnerEtAl_2011} --- ``since the CH is our ultimate desire, why don't we just shape it?'' --- the CG algorithm would have been a natural choice.

It is also evident that the sub-problem \eqref{SubProblem} is a profit maximization problem, 
given prices $\boldsymbol{\lambda}^{(k)}$ --- compare with \eqref{profit} and \eqref{Self-Scheduling}.
{\color{black}For clarity, sub-problem \eqref{SubProblem} is equivalently written as
\begin{equation} \label{maxprof}
\text{Sub}_i^{(k)}: (\mathbf{x}_{i}^S, \mathbf{y}_{i}^S; \boldsymbol{\lambda}^{(k)}) \in \underset{ (\mathbf{x}_{i}, \mathbf{y}_{i}) \in \mathcal{Z}_i}{\text{argmax}} \,\,  \phi_i (\mathbf{x}_i, \mathbf{y}_i; \boldsymbol{\lambda}^{(k)}),
\end{equation}
where $(\mathbf{x}_{i}^S, \mathbf{y}_{i}^S; \boldsymbol{\lambda}^{(k)})$ is the optimal self-schedule of unit $i$
given prices $\boldsymbol{\lambda}^{(k)}$.
Hence, using \eqref{SubProblem} and \eqref{maxprof}, 
the optimal profit under self-scheduling, $\phi_i (\mathbf{x}_{i}^S, \mathbf{y}_{i}^S; \boldsymbol{\lambda}^{(k)})$, is given by:
\begin{equation} \label{phi_h}
\phi_i (\mathbf{x}_{i}^S, \mathbf{y}_{i}^S; \boldsymbol{\lambda}^{(k)}) = - h_i^{(k)}(\mathbf{x}_{i}^S, \mathbf{y}_{i}^S).
\end{equation}
Using \eqref{phi_h} and \eqref{SubProblem}, the reduced cost for the optimal self-schedule $(\mathbf{x}_{i}^S, \mathbf{y}_{i}^S)$ from \eqref{rc} is:
\begin{align} \label{rc2}
rc_i^{(k)} (\mathbf{x}_{i}^S, \mathbf{y}_{i}^S) & = h_i^{(k)}(\mathbf{x}_{i}^S, \mathbf{y}_{i}^S) - \pi_i^{(k)} \nonumber \\
 & = - \pi_i^{(k)} - \phi_i (\mathbf{x}_{i}^S, \mathbf{y}_{i}^S; \boldsymbol{\lambda}^{(k)}). 
\end{align} }
The reduced cost calculation has a natural interpretation, in the spirit of \cite{BaumolFabian_1964},
which reveals a key property of CH prices and associated LOCs with respect to self-scheduling.
When unit $i$ ``is notified'' tentative prices $\boldsymbol{\lambda}^{(k)}$ at iteration $k$, 
the unit ``fictitiously'' self-schedules,
solving \eqref{SubProblem} --- equivalently \eqref{maxprof}  --- 
and calculating the tentative profit, 
$\phi_i (\mathbf{x}_i^S, \mathbf{y}_i^S; \boldsymbol{\lambda}^{(k)})$,
{\color{black}which equals $ - h_i^{(k)}(\mathbf{x}_{i}^S, \mathbf{y}_{i}^S)$. 
Standard duality of the LP problem in \eqref{UC_RMP} indicates that
the value of the dual $ - \pi_i^{(k)}$ represents the tentative profit of unit $i$, as calculated by the optimal selection of schedules in the RMP, 
which results in prices $\boldsymbol{\lambda}^{(k)}$. 
Given these prices $\boldsymbol{\lambda}^{(k)}$, 
if the unit can do better by self-scheduling, 
then, we have:
\begin{equation*}
\phi_i (\mathbf{x}_i^S, \mathbf{y}_i^S; \boldsymbol{\lambda}^{(k)}) > - \pi_i^{(k)}.  \end{equation*}
Hence, in this case, self-scheduling results in a negative reduced cost --- see \eqref{rc2}.}
The algorithm terminates when the profit-maximization problem \eqref{SubProblem} --- equivalently {\color{black}\eqref{maxprof}} --- yields a feasible schedule that already exists in the RMP for all units,
implying that for the CH prices, aggregate LOCs are minimized, and the optimal solution of the LD is reached.
{\color{black} Hence, upon termination, we derive the exact CH prices,
i.e., the prices that exactly solve the LD problem.}

\section{Illustration on Stylized Examples} \label{Stylized}

This section employs the stylized examples from \cite{SchiroEtAl_2016},
as a useful exercise, which aims at (i) illustrating the CG algorithm,
and (ii) providing intuition on the formation of CH prices viewed through the lens of D-W decomposition.
We refer the reader to \cite{SchiroEtAl_2016} for a detailed discussion of the CH price properties that each example illustrates.

\subsection{Example 1: 2-Gen, 1-Hour \cite[Ex. 1]{SchiroEtAl_2016}}
Two generators, G1 and G2, serve a 35 MW load, in a single period.
G1 is online, with technical minimum and maximum 10 and 50 MW, respectively, and energy offer \$50/MWh.
G2 can be either offline or online and submits a block offer of 50 MW at \$10/MWh. The MILP formulation is as follows:
\begin{subequations} \label{Ex1-milp}
\begin{equation} \label{Ex1-mip-obj}
\text{MILP: } \qquad \underset{x_1, x_2, y_2}{\min} \, f  =  50 x_1 + 10 x_2,  \qquad \qquad
\end{equation}
\begin{equation} \label{Ex1-mip1}
\text{subject to: } \qquad \qquad  x_1 + x_2 = 35, \qquad \qquad \qquad \quad
\end{equation}
\begin{equation} \label{Ex1-mip2}
10 \leq x_1 \leq 50,     
\end{equation}
\begin{equation} \label{Ex1-mip3}
    x_2 = 50 y_2, \quad y_2 \in \{0,1\}.
\end{equation}
\end{subequations}
This example has only one feasible solution with
G1 providing 35 MW, and MILP objective function value $f^* = \$1750$. The CH price is \$10/MWh, and G1 has \$1000 LOC.

Next, we solve this example using CG.
Consider trivial initial schedules: $z_1^{[1]}: \hat x_1^{[1]} = 10$, $\hat c_1^{[1]} = 500$,
and $z_2^{[1]}: \hat x_2^{[1]} = 0$, $\hat c_2^{[1]} = 0$.
{\color{black}To avoid RMP infeasibility, 
we introduce slack (deficit) variable $s$ to the RMP energy balance constraint \eqref{UCpbalRMP}, 
with a penalty in the objective function \eqref{UCobjRMP}.
Assuming a penalty value of \$1000/MWh,}
the initial RMP is as follows:
\begin{subequations} \label{Ex1-RMP}
\begin{align} 
\text{RMP}^{(1)}: \quad  \underset{z_1^{[1]}, z_2^{[1]}, s }{\min} g^{(1)} = 500 z_1^{[1]} + 0 z_2^{[1]} + 1000 s, \label{ex1RMPinit1}\\
\text{subject to: } \qquad 10 z_1^{[1]} + 0 z_2^{[1]} + s = 35, \rightarrow \lambda^{(1)}, \qquad \label{ex1RMPinit2} \\
z_1^{[1]} = 1, \rightarrow \pi_1^{(1)}, \qquad \qquad \qquad  \label{ex1RMPinit3a}\\ 
z_2^{[1]} = 1, \rightarrow \pi_2^{(1)}, \qquad \qquad \qquad \label{ex1RMPinit3b} 
\end{align}
\end{subequations}
with $z_1^{[1]}, z_2^{[1]}, s \geq 0$.
The solution of RMP$^{(1)}$ yields duals
$\lambda^{(1)} = 1000$, $\pi_1^{(1)} = -9500$, and $\pi_2^{(1)} = 0$. 
For these values, 
\begin{equation*} \label{ex1-G1-1}
\text{Sub}_1^{(1)}: \,\, \underset{10 \leq x_1 \leq 50}{\min} \, h_1^{(1)} = 50 x_1 - 1000 x_1 = -950 x_1,  \qquad \,\,
\end{equation*}
yields $x_1 = 50$, $h_1^{(1)} = -47500$, with $rc_1^{(1)} = h_1^{(1)} - \pi_1^{(1)} = -38000 < 0$, 
hence $z_1^{[2]}: \hat x_1^{[2]} = 50$, $\hat c_1^{[2]} = 2500$,
whereas
\begin{equation*} \label{ex1-G2-1}
\text{Sub}_2^{(1)}: \, \underset{x_2 = 50 y_2, y_2 \in \{0,1\}}{\min} \, h_2^{(1)} = 10 x_2 - 1000 x_2 = -990 x_2, \, 
\end{equation*}
yields $x_2 = 50$, $y_2 = 1$, $h_2^{(1)} = -49500$, with $rc_2^{(1)} = h_2^{(1)} - \pi_2^{(1)} = -49500 < 0$, 
hence $z_2^{[2]}: \hat x_2^{[2]} = 50$, $\hat c_2^{[2]} = 500$.
The RMP, after adding these two new columns becomes:
\begin{align*} \label{ex1RMPinit1-2}
\text{RMP}^{(2)}: \underset{\mathbf{z}, s}{\min} \, g^{(2)} & = 500 z_1^{[1]} { \color{black} + 0 z_2^{[1]} } \\
& + 2500 z_1^{[2]} + 500 z_2^{[2]} + 1000 s,
\end{align*}
\begin{equation*}\label{ex1RMPinit2-2}
\text {subject to: } \, 10 z_1^{[1]} { \color{black} + 0 z_2^{[1]} } + 50 z_1^{[2]} + 50 z_2^{[2]} + s = 35, \rightarrow \lambda^{(2)}, \quad
\end{equation*}
\begin{equation*}\label{ex1RMPinit3-2}
z_1^{[1]} + z_1^{[2]} = 1, \rightarrow \pi_1^{(2)},
\end{equation*}
\begin{equation*} \label{ex1RMPinit4-2}
z_2^{[1]} + z_2^{[2]} = 1, \rightarrow \pi_2^{(2)},
\end{equation*}
with $z_1^{[1]}, z_2^{[1]}, z_1^{[2]}, z_2^{[2]}, s \geq 0$.
The solution of RMP$^{(2)}$ now yields $\lambda^{(2)}$ = 10, $\pi_1^{(2)}$ = 400, and $\pi_2^{(2)}$ = 0. 
For these values, 
Sub$_1^{(2)}$ yields $x_1 = 10$, with $h_1^{(2)} = 400$, and $rc_1^{(2)} = h_1^{(2)} - \pi_1^{(2)} = 0$ 
--- note that a schedule with $x_1 = 10$ already exists in the RMP --- 
whereas Sub$_2^{(2)}$ yields $h_2^{(2)}  = 0$, for any feasible schedule, and $rc_2^{(2)}  = h_2^{(2)}  - \pi_2^{(2)} = 0$.
Since no feasible schedule with negative reduced cost is found, CG terminates.

The CH price is $\lambda^{(2)}$ = \$10/MWh.
The solution of RMP$^{(2)}$ is $z_1^{[1]} = 1$, $z_2^{[1]} = 0.5$, $z_2^{[2]} = 0.5$,
and the value of the objective function is $g^{(2)} = g^* = \$750$.
Hence, the duality gap between the MIP and the RMP is $f^* - g^* = \$1750 - \$750 = \$1000$, 
which represents the LOC to be paid to G1 (uplift).
Note that the maximum profit of G1 at the CH price would be derived by $x_1 = 10$, 
instead of $x_1 = 35$, which is the MILP dispatched quantity.
Notably, this quantity was the optimal solution for the G1 sub-problem when the algorithm terminated.
On the other hand, G2 has no LOC, since CH price equals its cost.

{\color{black}
\subsection{Example 2: Example 1 with 2 Nodes \cite[Ex. 2]{SchiroEtAl_2016}}

This example extends Example 1 to a two-node setting.
G1 is collocated with the 35 MW load at node 1, 
whereas G2 is located at node 2, 
and the transmission line between the two nodes has a capacity of 10 MW.
The MILP formulation follows \eqref{Ex1-milp} adding the transmission constraint: 
\begin{equation*}
    -10 \leq x_2 \leq 10.
\end{equation*}
The optimal solution remains the same (G1 provides 35 MW),
with no flow along the transmission line.
However, \cite[Ex. 2]{SchiroEtAl_2016} shows CH prices \$50/MWh at node 1, \$10/MWh at node 2, 
and a congestion price \$40/MWh.

Considering the initial schedules of the previous example,
and keeping, without loss of generality, 
only the one-direction flow that could be present in this example, 
i.e., $x_2 \leq 10$, 
the initial RMP follows \eqref{Ex1-RMP} 
adding the transmission constraint:
\begin{equation*} \label{ex1RMPinit5}
0 z_2^{[1]} \leq 10, \quad \longrightarrow \mu^{(1)}.
\end{equation*}
In fact, we should also include a slack variable, 
similarly to \eqref{ex1RMPinit2}, 
with an appropriate cost for the constraint violation in the objective function, 
but we omit this for the sake of brevity. 
Since the transmission constraint is not active in RMP$^{(1)}$, 
the same columns are generated and added to RMP$^{(2)}$, which now includes
the transmission constraint:
\begin{equation*} \label{ex1RMPinit5-2}
0 z_2^{[1]} + 50 z_2^{[2]} \leq 10, \quad \longrightarrow \mu^{(2)}.
\end{equation*}
The RMP$^{(2)}$ duals are $\lambda^{(2)} = 50$, $\mu^{(2)} = -40$, $\pi_1^{(2)} = 0$, and $\pi_2^{(2)} = 0$.
G1 sub-problem is the same with the previous example, 
with $rc_1^{(2)} = 0$,
whereas the objective function of G2 sub-problem, which now includes the transmission constraint dual, 
is given by $h_2^{(2)} =  10 x_2 - \lambda^{(2)} x_2 - \mu^{(2)} x_2  = 10 x_2 - 50 x_2 + 40 x_2 = 0$.     
Hence, $rc_2^{(2)} = h_2^{(2)} - \pi_2^{(2)} = 0$, and CG terminates 
with the aforementioned CH prices.
The solution of RMP$^{(2)}$ is $z_1^{[1]} =0.625$,  $z_1^{[2]} = 0.375$, $z_2^{[1]} = 0.8$, $z_2^{[2]} = 0.2$,
and the value of the objective function is $g^{(2)} = g^* = \$1350$.
Hence, the duality gap is $\$400$, 
representing the PRS due to the congestion price along the 10 MW transmission line.}

\subsection{Example {\color{black}3}: Example 1 with Start-up Cost \cite[Ex. 3]{SchiroEtAl_2016} }

This example adds to Example 1 a \$100 start-up cost to G2, 
i.e., $f_2(x_2, y_2) = 10 x_2 + 100 y_2$,
yielding a CH price equal to \$12/MWh that corresponds to the average cost of G2. 
Still, the only feasible solution is to dispatch G1 at 35 MW, 
and the MILP objective function value is $f^* = \$1750$.
Assuming same initial schedules, hence same RMP$^{(1)}$ duals, 
we get $z_1^{[2]}: \hat x_1^{[2]} = 50$, $ \hat c_1^{[2]} = 2500$, 
and $z_2^{[2]}: \hat x_2^{[2]} = 50$, $\hat c_2^{[2]} = 600$.
RMP$^{(2)}$ yields duals 
$\lambda^{(2)} = 12$, $\pi_1^{(2)} = 380$, and $\pi_2^{(2)} = 0$, 
for which no negative reduced cost schedule can be found. 
The solution of RMP$^{(2)}$ ($z$ variables) is the same with Example 1, 
but the value of the objective function is now $g^{(2)} = g^* = \$800$. 
Hence, the duality gap is $\$950$, also representing the G1 LOC.

\subsection{Example {\color{black}4}: 2-Gen, 2-Hour (linked) \cite[Ex. 4]{SchiroEtAl_2016}}

This example considers 2 hours, with load 45 and 80 MW.
G1 remains the same and should be online during both hours.
G2 has a 25 MW minimum, 
a 35 MW maximum, 
an energy offer of \$100/MWh, 
and should be either online or offline during both hours.
\begin{subequations} \label{Ex4-milp}
The MILP problem is formulated as follows:
\begin{equation} 
\text{MILP:} \,\underset{ \mathbf{x}, \mathbf{y} }{ \min } f = 50 x_{1,1} + 50 x_{1,2} + 100 x_{2,1} + 100 x_{2,2}, \label{Ex4-mip-obj}
\end{equation}
\begin{equation} 
 \text{subject to: } \quad   x_{1,1} + x_{2,1}  = 45, \quad   x_{1,2} + x_{2,2} = 80,  \label{Ex4-mip1-2}
 \end{equation}
\begin{equation}
\quad \quad \qquad 10 \leq x_{1,1} \leq 50,  \qquad 10 \leq x_{1,2} \leq 50,  \label{Ex4-mip2-2}
\end{equation}
\begin{equation}
25 y_2 \leq x_{2,1} \leq 35 y_2,   \quad 25 y_2 \leq x_{2,2} \leq 35 y_2, \label{Ex4-mip3}
\end{equation}
\begin{equation}
 y_2 \in \{0,1\},  \qquad \qquad  \qquad   \label{Ex4-mip4}
\end{equation}
\end{subequations}
yielding $x_{1,1} = 20$, $x_{1,2} = 50$, $x_{2,1} = 25$, $x_{2,2} = 30$, $y_2 = 1$, 
and MILP objective function value is $f^* = \$9000$.

Skipping the first CG steps, the last RMP is as follows:
\begin{align*}
\text{RMP}^{(3)}: \underset{ \mathbf{z}, \mathbf{s}}{\min} \,g^{(3)} = 1000 z_1^{[1]} + 0 z_2^{[1]}  + 5000 z_1^{[2]}  + 7000 z_2^{[2]} \,\, \\
 +  3000 z_1^{[3]} + 6000 z_2^{[3]} + 1000 s_1 + 1000 s_2, 
\end{align*}
subject to: 
\begin{align*}
10 z_1^{[1]} {\color{black}+ 0 z_2^{[1]}} + 50 z_1^{[2]} + 35 z_2^{[2]}  + 10 z_1^{[3]}  + & 25 z_2^{[3]} \\  
+ s_1 & = 45, \rightarrow \lambda_1^{(3)}, \\
10 z_1^{[1]}  {\color{black}+ 0 z_2^{[1]}} + 50 z_1^{[2]} + 35 z_2^{[2]}  + 50 z_1^{[3]}  + & 25 z_2^{[3]} \\ 
+ s_2 & = 80, \rightarrow  \lambda_2^{(3)}, \\
z_1^{[1]} + z_1^{[2]} + z_1^{[3]}  = 1, \rightarrow \pi_1^{(3)}, & \qquad \qquad \qquad \\
z_2^{[1]} + z_2^{[2]} + z_2^{[3]}  = 1, \rightarrow \pi_2^{(3)}, & \qquad \qquad \qquad
\end{align*}
with $z_1^{[1]}, z_2^{[1]}, z_1^{[2]}, z_2^{[2]}, z_1^{[3]}, z_2^{[3]}, s_1, s_2 \geq 0$.
RMP$^{(3)}$ yields duals $\lambda_1^{(3)} = 50$, $\lambda_2^{(3)} = 135.714$, $\pi_1^{(3)} = 4285.714$, and $\pi_2^{(3)} = 0$,
for which CG terminates.
The solution of RMP$^{(3)}$ is $z_1^{[2]} = 0.339$, $z_1^{[3]} = 0.661$, and $z_2^{[1]} = 0.143$, $z_2^{[3]} = 0.857$,
and the value of the objective function is $g^{(3)} = g^* = \$8821.429$.
Hence, the duality gap is $\$178.571$, 
representing G2 LOC.

\subsection{Example {\color{black}5}: Example {\color{black}4} with non-linked Hours \cite[Ex. 5]{SchiroEtAl_2016}}

This example considers the previous setting but as two single-period problems.
Equivalently, we can write the MILP formulation, 
using \eqref{Ex4-milp} and replacing \eqref{Ex4-mip3}--\eqref{Ex4-mip4} with:
$25 y_{2,1} \leq x_{2,1} \leq 35 y_{2,1}$,  
$25 y_{2,2} \leq x_{2,2} \leq 35 y_{2,2}$, 
and $y_{2,1}, y_{2,2} \in \{0,1\}$. 
CG terminates with RMP$^{(3)}$ yielding duals $\lambda_1^{(3)} = 50$, and  $\lambda_2^{(3)} = 100$, $\pi_1^{(3)} = -2500$, and $\pi_2^{(3)} = 0$, 
and a value of the objective function $g^{(3)} = g^* = \$7750$.
Hence the duality gap is $\$1250$, representing G2 LOC.
We leave this example as an exercise to the interested reader.

\section{{\color{black}Discussion and} More Realistic Test Cases} \label{MoreTests}

In this section, {\color{black}we discuss the comparative advantages of the proposed method and we proceed to more realistic test cases.
In Subsection \ref{Comp}, 
we present a comparative analysis with alternative approaches.
We support our analysis with two representative examples.
In Subsection \ref{Ramp}, we review an example provided in \cite{ChenEtAl_2020}, 
including ramp constraints, 
to illustrate the ease of accommodating features, 
which have been for long considered ``problematic.''
In Subsection \ref{RTO},
we test a 24-hour UC formulation on an ISO-sized FERC dataset \cite{FERC_RTO} with about 1000 generators, to illustrate scalability.}
 
{\color{black}
\subsection{Comparative Analysis of the D-W Characterization} \label{Comp}

The proposed method presents several advantages compared with the two main approaches in the computation of CH prices, 
namely sub-gradient methods and extended formulations.

The first approach, which refers to the initial attempts (e.g., \cite{WangEtAl_2009, WangEtAl_2013a}) to solve the Lagrangian Dual problem, was based on sub-gradient methods,
admittedly with not much success.
Arguably, sub-gradient methods are not well suited for the CH pricing problem,
due to the well-known convergence difficulties.
On the other hand, the D-W characterization combines several advantages.
It keeps the complexity in the small sub-problems 
--- as Lagrangian Relaxation also does --- 
but at the same time it shapes the CH, 
and comes with a convergence guarantee.
MILP-constrained feasible sets (with a polyhedral CH) guarantee finite convergence.

The second approach, which refers to more recent works 
(e.g., \cite{YuEtAl_2020, ChenEtAl_2020, KnuevenEtAl_2019}), 
has been directed to identifying extended formulations for the CH of each unit and the convex envelope of the cost function to derive an LP problem whose duals associated with the system constraints yield the CH prices.
Early works did not capture significant features, thus resulting in approximate descriptions.
For instance, ramping constraints were for long considered ``problematic,''
and could not be accommodated in these formulations.
Hence, significant effort was put on identifying extended formulations for ramp-constrained generators \cite{KnuevenEtAl_2018, BacciEtAl_2019} 
or applying disjunctive programming \cite{SchiroEtAl_2016} to cater for more general settings. 
As acknowledged in the relevant works, 
see e.g., \cite{YuEtAl_2020, KnuevenEtAl_2019}, 
extended formulations, 
even if they could accurately characterize specific classes of units,
they result in LP problems that are (at least) impractical to be solved in real-sized markets.
To address this issue, these works propose iterative algorithms, 
which start from formulations that describe exactly the CH of each unit, when possible, 
and approximate as tightly as possible the remaining ones.
The tightness of the approximations and their performance, 
i.e., the solution quality, and their speed are empirical questions, 
which depend on the employed formulation, the quality of their CH representation, etc.

The main characteristic of these latter works is that they rely on specific characterizations of special classes of units, and extensive reformulations,
on a case-by-case basis.
Notably, actual market implementations include far more complexity and details in the modeling and description of the unit characteristics and constraints.
Hence, such extensive reformulations impose a significant burden on their implementation in commercial products.
They also complicate modifications/additions of new units and/or constraints,
as they would require significant analysis to assess their impact on prices. 
Furthermore, extended formulations do not seem to provide an intuition of the CH price formation.

The proposed D-W characterization prevails as a natural fit to the CH pricing problem, 
and outperforms extended formulations in terms of simplicity, generalization, and economic interpretation.
It provides unique intuition about the CH price formation.
It does not require any reformulation of the unit constraints or cost functions.
These can be used ``as is,'' 
i.e., as modeled in current state-of-the-art Market Management Systems.
Even more, any change in the UC formulation, 
e.g., addition of new characteristics, new units, etc., 
can be integrated directly into our CH pricing problem formulation, unlike proposals in the literature to date that require case by case CH characterizations and/or approximations.
The proposed method is generalizable and independent of the commercial vendor-specific implementations of the unit costs/constraints.
This is an important feature and comparative advantage, 
which greatly facilitates regulatory actions 
and stakeholder acceptance.
Furthermore, it is a tractable exact method, with guaranteed finite convergence. 
Computational performance is always an empirical question.
The Operations Research theory based intuition is that the method is well-suited for the pricing problem
--- indeed the straightforward economic interpretation is not a coincidence.
Moreover, the computational experiments on an ISO-sized large-scale dataset confirms this intuition, 
and provides encouraging evidence of the scalability to actual market settings.

We further highlight some interesting computational features of the proposed method.
For the CG algorithm to proceed, 
we need to find at least one column with a negative reduced cost to update the duals.
That said, the sub-problems do not need to be solved to optimality in intermediate iterations, 
but only at the end, to guarantee termination.
Evidently, all sub-problems can be solved in parallel, 
or we may tailor the number of sub-problems solved at each iteration to the available computational resources.
Assuming full parallelization, 
the computational time of each iteration is dominated by the solution of the RMP, which is an LP problem,
including the system constraints and a ``convexity constraint'' for each unit.
The complexity of the unit specific constraints appears only in the sub-problems, not the RMP.
Using standard UC formulations, 
the sub-problem is a small MILP, usually very easy to solve.
Note that CG does not require to solve the profit maximization problem with a MILP formulation; 
we only need to find a negative reduced cost feasible schedule regardless of the method.
Nevertheless, a MILP formulation would definitely facilitate the implementation in current ISO Market Management Systems.}

\subsection{Ramp Constrained Example \cite[Ex. 2]{ChenEtAl_2020}} \label{Ramp}

This example includes 2 generators and a 3-hour horizon with load 95, 100, and 130 MW.
G1 has a maximum of 100 MW, and energy offer \$10/MWh. 
G2 has a 20 MW minimum, 35 MW maximum, energy offer \$50/MWh, start-up cost \$1000, no-load cost \$30, ramp rate 5 MW/hour, start-up rate 22.5 MW/hour, shut-down rate 35 MW/hour, minimum up/down times of 1 hour, and initially it is offline.
\begin{subequations} \label{3binRamp}
The UC 3-binary formulation is provided below \cite{ChenEtAl_2020}.
\begin{align}
     \underset{ \mathbf{p}, \mathbf{u}, \mathbf{v}, \mathbf{w} }{ \min } f = \sum_{t=1}^3 ( 10 p_{1,t}
    + 30 u_{2,t} + 50 p_{2,t} & + 1000 v_{2,t}),\\
\text{subject to: } \quad \qquad p_{1,t} + p_{2,t} = D_t, \quad & 1\leq t \leq 3,\\
0 \leq p_{1,t} \leq 100, \quad & 1\leq t \leq 3, \\
20 u_{2,t} \leq p_{2,t} \leq 35 u_{2,t}, \quad & 1\leq t \leq 3,\\
p_{2,t} - p_{2,t-1} \leq 5 u_{2,t-1} + 22.5 v_{2,t}, \quad & 1\leq t \leq 3,\\
p_{2,t-1} - p_{2,t} \leq 5 u_{2,t} + 35 w_{2,t}, \quad & 2\leq t \leq 3,\\
u_{2,t} - u_{2,t-1} = v_{2,t} - w_{2,t}, \quad & 1\leq t \leq 3,\\
v_{2,t} \leq u_{2,t}, \quad v_{2,t} \leq 1 - u_{2,t-1}, \quad & 1\leq t \leq 3,
\end{align}
\end{subequations}
with $p_{1,t}, p_{2,t} \geq 0$, and $u_{2,t}, v_{2,t}, w_{2,t} \in \{0,1\}$, $\forall t$, representing the status, start-up and shut-down variables, respectively.
The optimal schedule is:
$p_{1,1} = 75$, $p_{1,2} = 75$, $p_{1,3} = 100$, 
and $p_{2,1} = 20$, $p_{2,2} = 25$, $p_{2,3} = 30$, 
with $f^* = 7340$.

The interesting feature of this example is related to the ramping constraints,
which has been a challenge for the application of CH pricing.
The IR of \eqref{3binRamp} yields a vector of prices $\boldsymbol{\lambda}^\text{IR} = (10, 10, 182.701)$,
which deviate from the CH prices, which we show to be  $\boldsymbol{\lambda}^\text{CH} = (10, 10, 276)$.
We note that \cite{ChenEtAl_2020} reports Average Incremental Prices (AIC) prices 
obtained by an extended formulation that yields $\boldsymbol{\lambda}^\text{AIC} = (10, 10, 146.333)$,
while also noting that a numerical approximation using the 3-binary formulation 
would result in a very high price at hour 3, namely $(10, 10, 1161)$.
Apparently, in this example, hour 3 is of interest, and obviously, 
the IR (whose objective function value is $\$6464.55$) is less tight compared to the LD (whose value is shown to be $\$6975$).

Applying CG, {\color{black} using slack variables $s_t$, $\forall t$, with a penalty of \$1000/MWh}, we obtain: 
(1) Initial columns: $z_1^{[1]}: \mathbf{\hat p}_1^{[1]} = (0,0,0), \hat c_1^{[1]} = 0$; $z_2^{[1]}: \mathbf{\hat p}_2^{[1]} = (0,0,0), \hat c_2^{[1]} = 0$;
RMP$^{(1)}$ duals: $\boldsymbol{\lambda}^{(1)} = (1000, 1000, 1000)$, $\pi_1^{(1)}  = 0$, $\pi_2^{(1)}  = 0$.
(2) New columns: $z_1^{[2]}: \mathbf{\hat p}_1^{[2]} = (100, 100, 100)$, $\hat c_1^{[2]} = 3000$; $z_2^{[2]}: \mathbf{\hat p}_2^{[2]} = (22.5, 27.5, 32.5), \hat c_2^{[2]} = 5215$;
RMP$^{(2)}$  duals: $\boldsymbol{\lambda}^{(2)} = (-970, 0, 1000), \pi_1^{(2)} = 0, \pi_2^{(2)} = -5460$.
(3) New columns: $z_1^{[3]}: \mathbf{\hat p}_1^{[3]} = (0,0,100), \hat c_1^{[3]} = 1000; z_2^{[3]}: \mathbf{\hat p}_2^{[3]} = (0, 22.5, 27.5), \hat c_2^{[3]} = 3560$; 
RMP$^{(3)}$ duals: $\boldsymbol{\lambda}^{(3)} = (45.832, -25.832, 150.589), \pi_1^{(3)} = -14060, \pi_2^{(3)} = 0$.
(4) New columns: $z_1^{[4]}: \mathbf{\hat p}_1^{[4]} = (100,0,100), \hat c_1^{[4]} = 2000; 
z_2^{[4]}: \mathbf{\hat p}_2^{[4]} = (0, 0, 22.5), \hat c_2^{[4]} = 2155$.
RMP$^{(4)}$ duals: $\boldsymbol{\lambda}^{(4)} = (10, 10, 276), \pi_1^{(4)} = -26600, \pi_2^{(4)} = -4255$.
CG terminates with $z_1^{[2]} = 0.75, z_1^{[3]}=0.163, z_1^{[4]} = 0.088$, and $z_2^{[2]} = 0.5, z_2^{[3]} = 0.5$; objective function value:  $g^{(4)} = g^* = \$6975$.

\begin{figure}[tb]
\centering
\includegraphics[width=3.45in]{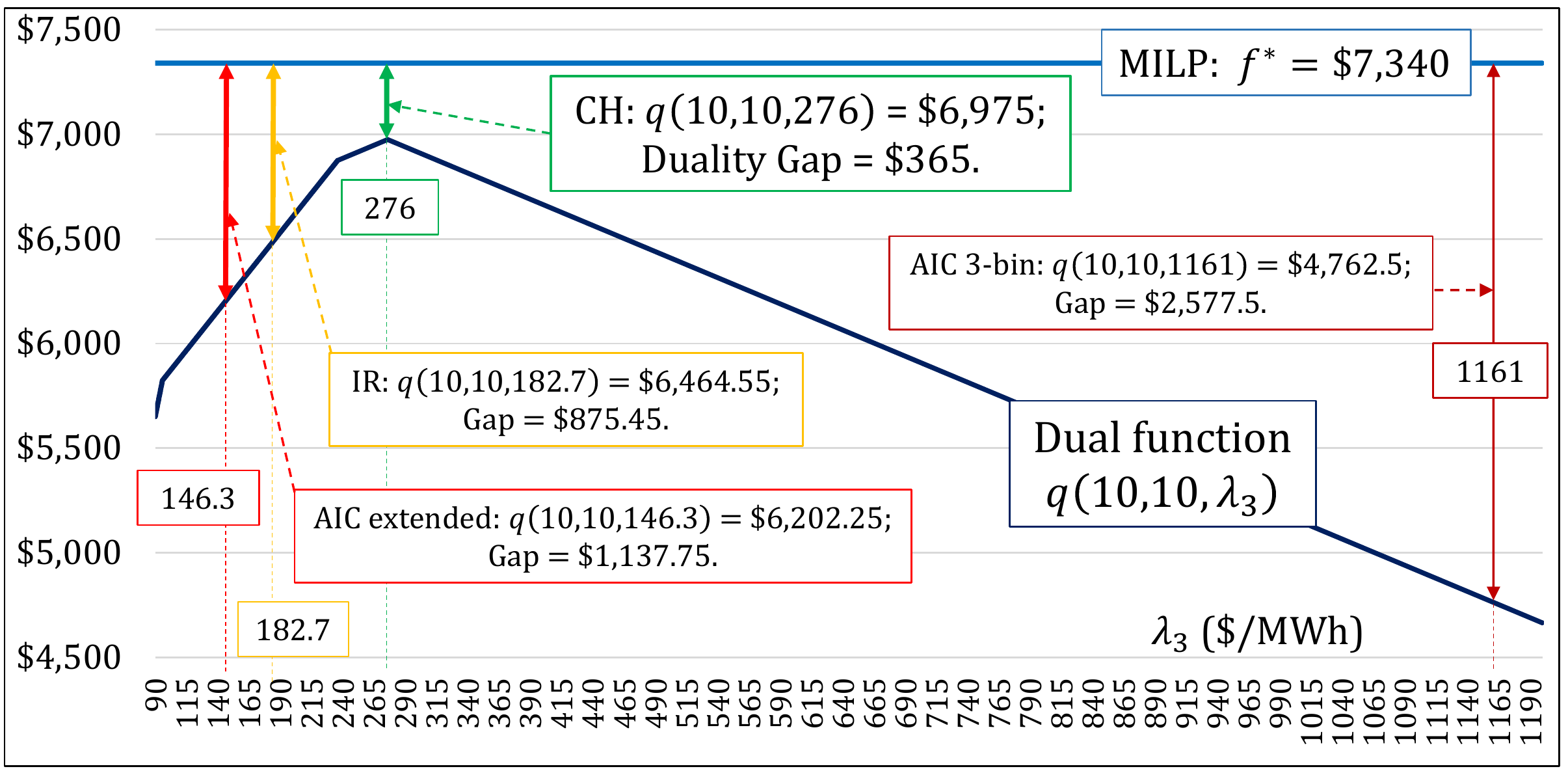}
\caption{Evaluation of $q(\boldsymbol{\lambda})$, for $\boldsymbol{\lambda} = (10, 10, \lambda_3)$, $90 \leq \lambda_3 \leq 1200.$} 
\label{Fig1}
\end{figure}

In Fig. \ref{Fig1}, we illustrate the evaluation of the dual function $q(\boldsymbol{\lambda})$ for $\boldsymbol{\lambda} = (10, 10, \lambda_3)$,
with $\lambda_3$ (\$/MWh) ranging from 90  --- the marginal cost reported in \cite{ChenEtAl_2020} --- to 1200, so as to include the high price of 1161.
The CH duality gap is $\$365$, and corresponds to LOC of G2.
Notably, the IR gap is $\$875.45$,
whereas the extended formulation AIC gap is $\$1137.75$.

{\color{black}Lastly, we also implemented a sub-gradient method, 
which solves the same sub-problems with CG, 
but updates the duals (prices $\boldsymbol{\lambda}$) using a step-size.
Not surprisingly, step-size tuning is not a trivial task.
We experimented with several step-sizes, 
namely $1/k, 1/\sqrt{k}, 10/k, 10/\sqrt{k}, 0.001, 0.01, 0.02$,
and different initializations, 
e.g., $\boldsymbol{\lambda}^{(1)} = (10, 10, 90)$ (the marginal cost),
and $\boldsymbol{\lambda}^{(1)} = \boldsymbol{\lambda}^\text{IR} (10, 10, 182.701)$ (the IR prices).
Many of the step-sizes exhibited poor performance. 
Our best results were obtained for the step-size $10/k$, 
but, even for this small example, 
we only managed to reach prices to a proximity of about 0.1 (\$/MWh) for all three hours, 
after 10,000 iterations.
For comparison, a typical $1/k$ step-size yielded (after 10,000 iterations)
a price for $\lambda_3$ of about 206 (\$/MWh),
i.e., 70 (\$/MWh) less than the exact CH price.
In brief, CG significantly outperformed the sub-gradient implementation in the following aspects:
\emph{(i)} CG provides the exact LD values, whereas a sub-gradient method provides approximate bounds; 
\emph{(ii)} CG converges finitely, 
whereas a sub-gradient method may only converge asymptotically 
after orders of magnitude more iterations, and 
\emph{(iii)} CG provides exact prices without need of tuning, 
whereas a sub-gradient method provides only approximate prices if properly tuned, 
which as we show is a non-trivial task even for a small example.
}

\subsection{FERC Dataset Example} \label{RTO}

In this subsection, 
we test our approach on the PJM-like FERC ``summer'' dataset \cite{FERC_RTO}.
UC formulations have been extensively studied 
--- see e.g., \cite{KnuevenEtAl_2020} for an overview of MILP formulations,
and our goal in this last example is not to exhaustively include all features,
but to provide evidence on the scalability to real-sized problems.
In fact, as it has already become evident,
a major strength of CG is that it can easily accommodate any unit model, 
and all system constraints currently present in electricity markets.

\begin{subequations} \label{UCform}
The employed UC MILP formulation is provided below.
Unless, otherwise mentioned, $i \in \mathcal{I}$, $t\in \mathcal{T}$, and $b \in \{1,...B\}$,
where $B$ is the number of block offer steps.

\begin{align}
\underset{ \mathbf{p}, \mathbf{u}, \mathbf{v}, \mathbf{w}, \mathbf{s}}{ \min } f = 
    \sum_{i, t} \Big[ 
    C_{i}^{NL}u_{i,t} 
    + C_{i}^{SU} v_{i,t} 
    + C_{i}^{SD} w_{i,t} \nonumber \\
    + \sum_b \big( C_{i,b}^p \, p_{i,t,b} \big)
    + C_{i}^r \, r_{i,t} \Big]
    + \sum_{t} \big(
     M^p s_t^{p}
    + M^r s_t^{r} \big), \label{UCform-obj}\\
\text{subject to: } \qquad \sum_{i} p_{i,t} + s_t^{p} = D_t^p, \quad \forall t, \qquad \label{EnBal} \\
\sum_{i} r_{i,t} + s_t^{r} \geq D_t^r, \quad \forall t, \qquad \label{ResReq} \\
p_{i,t,b} \leq {\bar P}_{i,b} \, u_{i,t}, \, \, \forall i,t,b, \quad p_{i,t} = \sum_{b} p_{i,t,b}, \quad \forall i,t, \label{Block}\\
{\underline P}_i \, u_{i,t} \leq p_{i,t} \leq {\bar P}_i \, u_{i,t} - r_{i,t} , \quad r_{i,t} \leq {\bar R}_{i}, \quad \forall i,t,  \label{Capacity}\\
p_{i,t} - p_{i,t-1} \leq R_i^U u_{i,t-1} + R_i^{SU} v_{i,t}, \quad \forall i, t, \label{RampU}\\
p_{i,t-1} - p_{i,t} \leq R_i^D u_{i,t} + R_i^{SD} w_{i,t}, \quad \forall i, t, \label{RampD}\\
u_{i,t} - u_{i,t-1} = v_{i,t} - w_{i,t}, \quad \forall i, t, \label{SUSD} \\
\sum_{t' = t - {MUT}_i + 1}^t v_{i,t'} \leq u_{i,t},  \quad \forall i, t = [MUT_i, T],\label{MU} \\
\sum_{t' = t - {MDT}_i + 1}^t w_{i,t'} \leq 1 - u_{i,t},  \quad \forall i, t = [MDT_i, T], \label{MD}\\
\sum_{t' = 1}^{UT_i} u_{i,t'} = UT_i, \quad \sum_{t' = 1}^{DT_i} u_{i,t'} = 0, \quad \forall i, \label{UTDT}
\end{align}
with $p_{i,t, b}, p_{i,t}, r_{i,t}, s_t^{p}, s_t^{r} \geq 0$, and $u_{i,t}, v_{i,t}, w_{i,t} \in \{0,1\}$, $\forall i,t,b$.

In brief, \eqref{UCform-obj} minimizes the aggregate commitment costs (no-load, $C_i^{NL}$, start-up, $C_i^{SU}$, and shut-down, $C_i^{SD})$, dispatch costs ($C_{i,b}^p$ is the cost of block offer $b$),\footnote{{\color{black} Please note that the block offers (with an increasing cost $C_{i,b}^p$, as per market rule requirements) yield a piecewise linear objective function.}} 
plus the reserve cost (reserve offer $C_i^r$),
plus a penalty cost related to power/reserve deficit slacks $s_t^p$, $s_t^r$, with penalties $M^p$ and $M^r$.
System constraints \eqref{EnBal}--\eqref{ResReq} include power balance (demand $D_t^p$) and reserve requirements ($D_t^r$).
Unit specific constraint sets \eqref{Block}--\eqref{UTDT} are as follows.
\eqref{Block} imposes block offer $b$ maximum, ${\bar P}_{i,b}$, and defines $p_{i,t}$ as the sum of accepted block quantities $p_{i,t,b}$.
\eqref{Capacity} imposes minimum and maximum capacity limits for power, ${\underline P}_i$, ${\bar P}_i$, and maximum reserve capability ${\bar R}_i$.
\eqref{RampU}--\eqref{RampD} impose ramp up/down limits $R_i^U$, $R_i^D$, including start-up/shut-down, $R_i^{SU}$, $R_i^{SD}$.
\eqref{SUSD} defines commitment variables, 
\eqref{MU}--\eqref{MD} impose the minimum up/down time limits, $MUT_i$, $MDT_i$, 
whereas \eqref{UTDT} enforces initial up/down times; $UT_i$ and $DT_i$ are the number of hours the unit should remain online or offline, respectively, depending on initial conditions
(truncated at $T$).
\end{subequations}

The RMP at iteration $k$, RMP$^{(k)}$, is formulated as:
\begin{align*}
\underset{\mathbf{z}, \mathbf{s}}{\min} \,g^{(k)} =  \sum_{i, n_i \in \mathcal{N}_i^{(k)}} \hat c_i^{[n_i]}z_{i}^{[n_i]} 
    & + \sum_{t} ( M^p s_t^{p} + M^r s_t^{r} ), 
\end{align*}
\begin{align*}
\text{subject to:} \sum_{i, n_i \in \mathcal{N}_i^{(k)}} \hat p_{i,t}^{[n_i]}z_{i}^{[n_i]} + s_t^{p} & = D_t^p, \quad \forall t \rightarrow \lambda_t^{p\,(k)}, \\
\sum_{i, n_i \in \mathcal{N}_i^{(k)}} \hat r_{i,t}^{[n_i]}z_{i}^{[n_i]}  + s_t^{r} & \geq D_t^r, \quad \forall t, \rightarrow \lambda_t^{r\,(k)},\\
\sum_{n_i \in \mathcal{N}_i^{(k)}} z_{i}^{[n_i]}  & = 1,\qquad  \forall i \rightarrow \pi_i^{(k)},
\end{align*} 
with $z_{i}^{[n_i]} \geq 0, \, \forall i, \, n_i \in \mathcal{N}_i^{(k)}$, and $s_t^p, s_t^r \geq 0, \, \forall t  $.

The unit $i$ sub-problem at iteration $k$, Sub$_i^{(k)}$, is given by:
\begin{align*}
 \min \, h_i^{(k)} =    \sum_{t} \Big[
    C_{i}^{NL}u_{i,t} 
    + C_{i}^{SU} v_{i,t} 
    + C_{i}^{SD} w_{i,t}   \qquad \\
     + \sum_{b} \big( C_{i,b}^p \, p_{i,t,b} \big) + C_i^r \, r_{i,t} \Big]
     - \sum_{t} 
    \big(
           \lambda_t^{p \, (k)}  p_{i,t}
         + \lambda_t^{r \, (k)}  r_{i,t}  
    \big), 
\end{align*}
subject to: unit specific constraints \eqref{Block}--\eqref{UTDT}. 
The reduced cost is $ rc_i^{(k)}  =  h_i^{(k)}  - \pi_i^{(k)}$.

The dataset was adjusted to fit \eqref{UCform}.
It contains 1011 units, {\color{black}each with up to 10-block energy offers} and zero-priced reserve offers.
The 24-hour $D_t^p$ was set equal to forecasted demand, 
and penalties $M^p$ and $M^r$ to 1000 and 900.
Hot-type values were used for $C_i^{SU}$ ,
missing $MUT_i$ and $MDT_i$ were set to 1,
$UT_i$ and $DT_i$ to 0, 
and all units initially online dispatched at $\underline{P}_i$.
$R_i^U$ and $R_i^D$ were reduced to half, 
$R_i^{SU}$ were set equal to $\underline{P}_i$ plus half the $R_i^U$, 
$R_i^{SD}$ equal to $\bar{P}_i$,
and $\bar{R}_i$ equal to 30-min ramp up.
CG was modeled in C, and solved on a Intel Core i7-5500U \@2.4GHz with 8 GB RAM, using CPLEX 12.7.

\begin{figure}[tb]
\centering
\includegraphics[width=3.45in]{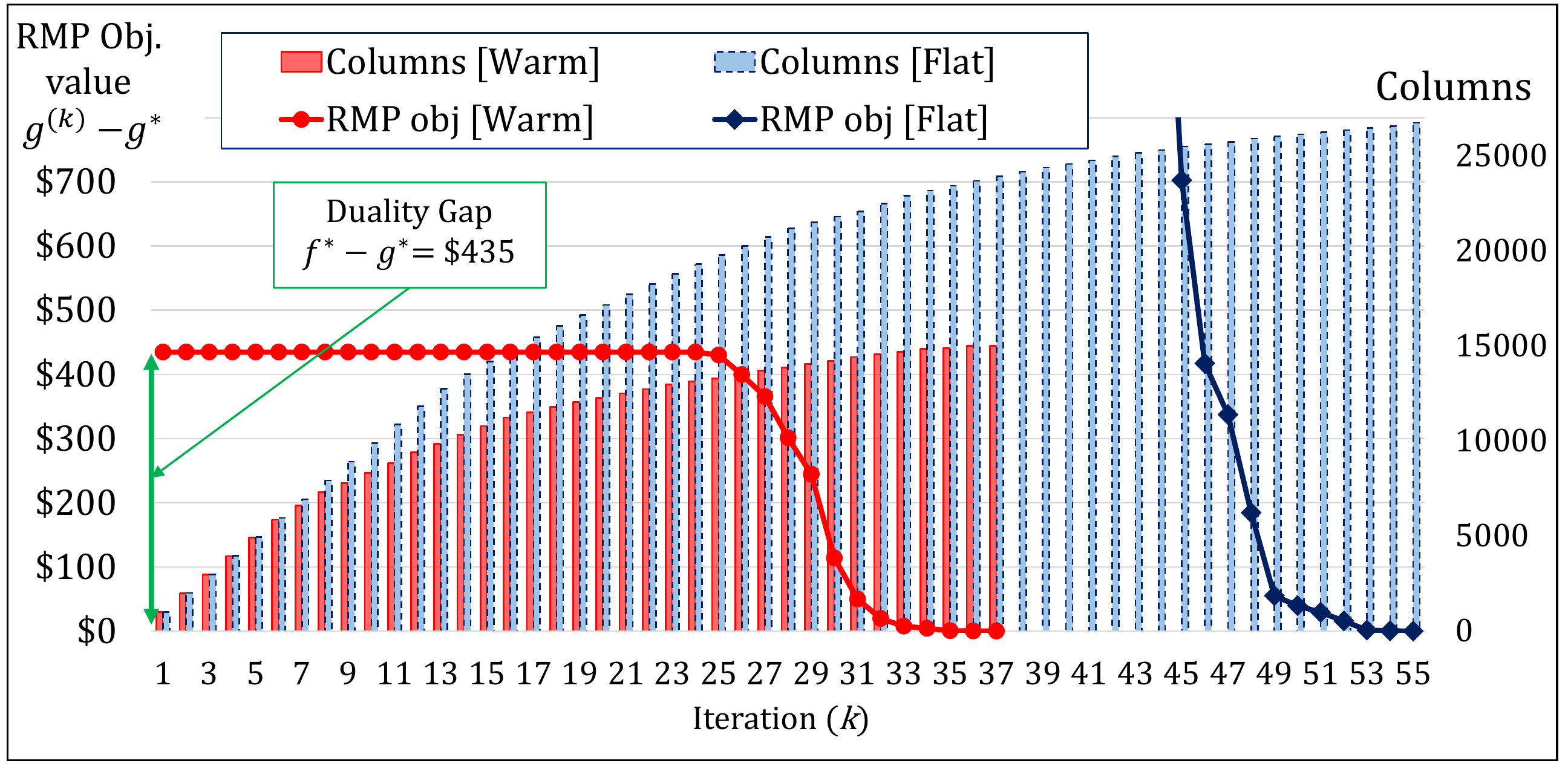}
\caption{CG iterations for flat and warm start. Left axis shows the RMP objective function value as the difference of $g^{(k)}$ from the optimal LD bound $g^* = \$32,169,154.7$. Right axis shows total columns at each iteration. UC MILP solved with 0.01\% optimality gap. For comparison purposes, note that the objective function value of IR is $\$32,149,006.2$ (less tight bound).}
\label{Fig2}
\end{figure}
In Fig. \ref{Fig2}, we explore CG performance in terms of iterations, 
when initiated from (i) the MILP UC solution --- referred to as a ``Warm start,'' 
and (ii) a ``Flat start'' with initial columns assuming self-scheduling at zero prices.
Note that current practice normally performs the pricing run after the MILP UC solution.
However, CH pricing could be run in parallel to the MILP UC problem,
as it does not actually need its solution! 
The results indicate that CG terminates in a few tens of iterations for both cases.
We do observe an anticipated plateau at the beginning of the warm start, 
since CG initially needs some columns to be able to ``combine'' them.
The flat start drops rapidly (high values due to deficit variable penalties are truncated), 
but it takes a few more iterations to terminate.
The computational times in a sequential implementation are dominated by the sub-problems; 
each iteration took about 14 seconds.
However, we note that the RMP time was less than 1 second (usually significantly less),
whereas each sub-problem was solved in the order of 10 msec.
It goes without saying that parallel implementation of sub-problems would drive down each iteration to practically the time to solve the RMP,
i.e., an LP problem; the estimated time for both cases assuming full parallelization would be well below 1 minute.
But even with some simple partial parallelization, say e.g., solving a batch of a dozen of sub-problems,
the time per iteration would be reduced by approximately a factor of 10.

\section{Conclusions and Further Research} \label{Conclusions}

In this work, we compute CH prices, employing D-W decomposition and CG. 
We balance the narrative with describing the theory,
enhancing the reader's intuition with illustrative examples, 
and providing indications of computational tractability and scalability to real-sized systems.
The simplicity of CG,
the intuitive explanation of the problem ``convexification'' or CH formation,
and the amenable to parallelization structure using standard UC and unit specific MILP formulations
reinforce the use of CG for explanation and computation.
There was no use of any of the known methods that may enhance CG performance, although such strategies have been successful in other applications. Ongoing work is directed in testing CG with actual market data on a large US ISO application.


\ifCLASSOPTIONcaptionsoff 
  \newpage
\fi


\begin{thebibliography}{99}

\bibitem{CaramanisEtAl_1982}
M.~C.~Caramanis, R.~R.~Bohn, F.~C.~Schweppe. 
``Optimal spot pricing: Practice and theory,'' 
\emph{IEEE Trans. Power App. Syst.},
vol. PAS-101, no. 9, pp. 3234--3245, 1982.

\bibitem{SchweppeEtAl_1988}
F.~C.~Schweppe, M.~C.~Caramanis, R.~D.~Tabors, R.~E.~Bohn,
``\emph{Spot pricing of electricity},'' Kluwer Academic Publishers, Boston, MA, 1988.

\bibitem{ONeillEtAl_2005}
R.~P.~O'Neill, P.M.~Sotkiewicz, B.~F.~Hobbs, M.~H.~Rothkopf, and W.~R.~Stewart~Jr.,
``Efficient market-clearing prices in markets with nonconvexities,''
\emph{Europ. J. Oper. Res.}, vol 164, no. 1, pp. 269--285, 2005.

\bibitem{AndrianesisEtAl_2013a}
P.~Andrianesis, G.~Liberopoulos, G.~Kozanidis, and A.~D.~Papalexopoulos,
``Recovery mechanisms in day-ahead markets with non-convexities --- Part I: Design and evaluation methodology,''
\emph{IEEE Trans. Power Syst.}, vol 28, no. 2, pp. 960--968, 2013.

\bibitem{AndrianesisEtAl_2013b}
P.~Andrianesis, G.~Liberopoulos, G.~Kozanidis, and A.~D.~Papalexopoulos,
``Recovery mechanisms in day-ahead markets with non-convexities --- Part II: Implementation and numerical evaluation,''
\emph{IEEE Trans. Power Syst.}, vol 28, no. 2, pp. 969--977, 2013.

\bibitem{HoganRing_2003}
W.~W.~Hogan and B.~J.~Ring, ``On minimum-uplift pricing for electricity markets,'' Working Paper, John F. Kennedy School of Government, Harvard University, 2003.

\bibitem{GribikEtAl_2007}
P.~R.~Gribik, W.~W.~Hogan, and S.~L.~Pope, ``Market-clearing electricity prices and energy uplift,'' Working Paper, John F. Kennedy School of Government, Harvard University, 2007.

\bibitem{ONeillEtAl_2017}
R.~P.~O'Neill, A.~Castillo, B.~Eldridge, and R.~B.~Hytowitz,
``Dual algorithm in ISO markets,''
\emph{IEEE Trans. Power Syst.}, vol 32, no. 4, pp. 3308--3310, 2017.

\bibitem{LibAnd_2016}
G.~Liberopoulos and P.~Andrianesis, ``Critical review of pricing schemes in markets with non-convex costs,'' 
\emph{Oper. Res.}, vol. 64, no. 1, pp. 17-31, 2016.

\bibitem{FERC_2014}
Federal Energy Regulatory Commission, 
``Price formation in organized wholesale electricity markets,''
Docket No. AD14-14-000, Dec. 2014.

\bibitem{WangEtAl_2016}
C.~Wang, P.~B.~Luh, P.~Gribik, T.~Peng, and L.~Zhang, 
``Commitment cost allocation of fast-start units for
approximate extended locational marginal prices,''
\emph{IEEE Trans. Power Syst.}, vol 31, no. 6, pp. 4176--4184, 2016.

\bibitem{ChenWang_2017}
Y.~Chen, and C.~Wang, 
``Enhancements of extended locational marginal pricing --- Advancing practical implementation,''
Midcontinent Independent System Operator, Nov. 2017.

\bibitem{PJM_2018}
PJM, 
``Important concepts from price formation education session 2: Alternative
pricing frameworks,'' 2018.

\bibitem{MISO_2019}
Midcontinent Independent System Operator, 
``ELMP III white paper I, R\&D report and design recommendation on short-term enhancements,'' Jan. 31, 2019.

\bibitem{Chao_2019}
H.~Chao,
``Incentive for efficient pricing mechanism in markets with non-convexities,''
\emph{J. Reg. Econ.}, vol 56, pp. 33--58, 2019.

\bibitem{ONeill_2020}
R.~P.~O'Neill,
``Nonconvex electric power auction markets: market efficiency and pricing,'' 
presented at INFORMS Annual Meeting 2020.

{\color{black}
\bibitem{bertsim}
D.~Bertsimas, E.~Litvinov, X.~A.~Sun, J.~Zhao and T.~Zheng, 
``Adaptive robust optimization for the security constrained unit commitment problem,'' 
\emph{IEEE Trans. Power Syst.}, vol. 31, no. 6, pp. 52--63, 2013.

\bibitem{ding1}
T. Ding, Z. Wu, J. Lv, Z. Bie and X. Zhang, 
``Robust co-optimization to energy and ancillary service joint dispatch considering wind power uncertainties in real-time electricity markets,'' 
\emph{IEEE Trans. Sustain. Energy}, vol. 7, no. 4, pp. 1547--1557, 2016.

\bibitem{ding2}
T. Ding et al., 
``Duality-free decomposition based data-driven stochastic security-constrained unit commitment,''
\emph{IEEE Trans. Sustain. Energy}, vol. 10, no. 1, pp. 82--93, 2019.
}
\bibitem{WangEtAl_2009}
C.~Wang, P.~B.~Luh, P.~Gribik, L.~Zhang, and T.~Peng,
``A subgradient-based cutting plane method to calculate convex hull market prices,''
in \emph{Proc. 2009 IEEE PES GM,}, Calgary, AB, Canada, 26--30 July 2009.

\bibitem{WangEtAl_2013}
C.~Wang, T.~Peng, P.~B.~Luh, P.~Gribik, and L.~Zhang,
``The subgradient Simplex cutting plane method for extended locational marginal prices,''
in \emph{IEEE Trans. Power Syst.}, vol. 28, no. 3, pp. 2758--2767, 2013.

\bibitem{WangEtAl_2013a}
G.~Wang, U.~V.~Shanbhag, T.~Zheng, E.~Litvinov, and S.~Meyn,
``An extreme-point subdifferential method for convex hull pricing in energy and reserve markets --- Part I: Algorithm structure,''
\emph{IEEE Trans. Power Syst.}, vol 28, no. 3, pp. 2111--2120, 2013.

\bibitem{WangEtAl_2013b}
G.~Wang, U.~V.~Shanbhag, T.~Zheng, E.~Litvinov, and S.~Meyn,
``An extreme-point subdifferential method for convex hull pricing in energy and reserve markets --- Part II: Convergence analysis and numerical performance,''
\emph{IEEE Trans. Power Syst.}, vol 28, no. 3, pp. 2121--2127, 2013.

\bibitem{ItoEtAl_2013}
N.~Ito, A.~Takeda and T.~Namerikawa,
``Convex hull pricing for demand response in electricity markets,''
in \emph{Proc. IEEE SmartGridComm}, Vancouver, BC, 21-24 Oct. 2013.

\bibitem{HuaBaldick_2017}
B.~Hua and R.~Baldick,
``A convex primal formulation for convex hull pricing,''
\emph{IEEE Trans. Power Syst.}, vol 32, no. 5, pp. 3814--3823, 2017.

\bibitem{GarciaEtAl_2020}
M.~Garcia, H.~Nagarajan, and R.~Baldick,
``Generalized convex hull pricing for the {AC} optimal power flow problem,''
\emph{IEEE Trans. Control Netw. Syst.}, vol. 7, no. 3, pp. 1500--1510, 2020.

\bibitem{YangEtAl_2019}
Z.~Yang, T.~Zheng, J.~Yu, and K.~Xie,
``A unified approach to pricing under nonconvexity,''
\emph{IEEE Trans. Power Syst.}, vol 34, no. 5, pp. 3417--3427, 2019.

\bibitem{YuEtAl_2020}
Y.~Yu, Y.~Guan, and Y.~Chen,
``An extended integral unit commitment formulation and an iterative algorithm for convex hull pricing,''
\emph{IEEE Trans. Power Syst.}, vol 35, no. 6, pp. 4335--4346, 2020.

\bibitem{ChenEtAl_2020}
Y.~Chen, R.~O'Neill, and P.~Whitman, ``A Unified approach to solve convex hull pricing and average incremental cost pricing with large system study,'' Working Paper, 2020.

\bibitem{AlvarezEtAl_2020}
C.~Alvarez, F.~Mancilla-David, P.~Escalona, and A.~Angulo,
``A Bienstock-Zuckerberg-based algorithm for solving a network-flow formulation of the convex hull pricing problem,''
\emph{IEEE Trans. Power Syst.}, vol 35, no. 3, pp. 2108--2114, 2020.

\bibitem{KnuevenEtAl_2019}
B.~Knueven, J.~Ostrowski, A.~Castillo, and J.-P.~Watson,
``A computationally efficient algorithm for computing convex hull prices,'' SAND2019-10896 J, Sandia National Laboratories, Albuquerque, NM, Sep. 2019.

\bibitem{SchiroEtAl_2016}
D.~A.~Schiro, T.~Zheng, F.~Zhao, and E.~Litvinov, ``Convex hull pricing in electricity markets: Formulation, analysis, and implementation challenges,'' \emph{IEEE Trans. Power Syst.}, vol. 31, no. 5, pp. 4068--4075, 2016.

\bibitem{AndKoz_2014}
P.~Andrianesis and G.~Kozanidis, 
``A multi-stage column generation solution approach for the bidline aircrew scheduling problem,''
in Proc. 3rd Intl. Symp. Oper. Res., Volos, Greece June 26--28, 2014.

\bibitem{DantzigWolfe_1960}
G.~B.~Dantzig and P.~Wolfe, ``Decomposition Principle for Linear Programs,'' \emph{Oper. Res.}, vol. 8, no. 1, pp. 101-111, 1960.

\bibitem{Geoffrion_1974}
A.~M.~Geoffrion, ``Lagrangian relaxation for integer programming,'' 
\emph{Mathem. Program. Study}, pp. 82--114, 1974.

\bibitem{MagnantiEtAl_1976}
T.~L.~Magnanti, J.~F.~Shapiro, and M.~H.~Wagner, 
``Generalized linear programming solves the dual,''
\emph{Manag. Sci.}, vol. 22, no. 11, pp. 1195--1203, 1976.

\bibitem{ManneEtAl_1980}
A.~Manne, H.~Chao, and R.~Wilson,
``Computation of competitive equilibria by a sequence of linear programs,''
\emph{Econometrica}, vol. 48, no. 7, pp. 1595--1615, 1980.

\bibitem{BarnhartEtAl_1998}
C.~Barnhart, E.~L.~Johnson, G.~L.~Nemhauser, M.~W.~P.~Savelsebergh, and P.~H.~Vance, 
``Branch-and-Price: Column generation for solving huge integer programs,''
\emph{Oper. Res.}, vol. 46, no. 3, pp. 316--329, 1998.

\bibitem{Vanderbeck_2000}
F.~Vanderbeck, 
``On Dantzig-Wolfe decomposition in integer programming and ways to perform branching in a branch-and-price algorithm,''
\emph{Oper. Res.}, vol. 48, no. 1, pp. 111--128, 2000.

\bibitem{LubbeckeDesrosiers_2005}
M.~E.~L{\"u}bbecke and J.~Desrosiers,
``Selected Topics in Column Generation,'' 
\emph{Oper. Res.}, vol. 53, no. 6, pp. 1007--1023, 2005.

\bibitem{FERC_RTO}
E.~Krall, M.~Higgins, and R.~P.~O'Neill, 
``RTO unit commitment test system,''
FERC Staff Report, 2012.

\bibitem{FeltKiwiel_2001}
S.~Feltenmark and K.~C.~Kiwiel,
``Dual applications of proximal bundle methods, including Lagrangian relaxation of nonconvex problems,''
\emph{SIAM J. Optim.}, vol. 10, no. 3, pp. 697--721, 2001.

\bibitem{Lemarechal_2007}
C.~L{\'e}marechal,
``The omnipresence of Lagrange,'' 
\emph{Ann. Oper. Res.}, vol. 153, pp. 9--27, 2007.

\bibitem{LunaEtAl_2020}
J.~P.~Luna, C.~Sagastiz\'abal, and P.~J.~S.~Silva,
``A discussion on electricity prices, or the two sides of the coin,'' 2020,
http://www.optimization-online.org/DB\_HTML/2020/07/7904.html.


\bibitem{BertsekasEtAl_1983}
D.~P. Bertsekas, G.~S.~Lauer, N.~R.~Sandell, and T.~A.~Posberg,
``Optimal short-term scheduling of large-scale power systems,'' 
\emph{IEEE Trans. Aut. Control}, vol. 28, pp. 1--11, 1983.

\bibitem{Dupin_2019}
N.~Dupin,
``Column generation for the discrete {UC} problem with min-stop ramping constraints,''
\emph{IFAC PapersOnLine}, vol 52, no. 13, pp. 529--534, 2019.

\bibitem{Ber}
D. P. Bertsekas, 2016. Nonlinear Programming. Third Edition. Athena Scientific, Belmont, MA.

\bibitem{KnuevenEtAl_2018}
B.~Knueven, J.~Ostrowski, and J.~Wang, ``The ramping polytope and cut generation for the unit commitment problem,'' \emph{INFORMS J. Comput.}, vol. 30, no. 4, pp. 739--749, 2018.

\bibitem{BacciEtAl_2019}
T.~Bacci, A.~Frangioni, C.~Gentile, and K.~Tavlaridis-Gyparakis, 
``New MINLP formulations for the unit commitment problems with ramping constraints,'' 2019, http://www.optimization-online.org/DB\_HTML/2019/10/7426.html.

\bibitem{Falk_1969}
J.~E.~Falk, ``Lagrange multipliers and nonconvex programs.'' \emph{SIAM J. Control}, vol. 7, no. 4, pp. 534–-545, 1969.

\bibitem{BergnerEtAl_2011}
M.~Bergner, A.~Caprara, F.~Furini, M.~E.~L{\"u}bbecke, E.~Malaguti, and E.~Traversi.
``\emph{Partial Convexification of General MIPs by Dantzig-Wolfe Reformulation},'' 
in O.~G{\"u}nl{\"u}k, G.~J.~Woeginger (eds) Integer Programming and Combinatoral Optimization.
Lecture Notes in Computer Science, vol. 6655. Springer, Berlin, Heidelberg, 2011.

\bibitem{BaumolFabian_1964}
W.~J.~Baumol and T.~Fabian,
``Decomposition, pricing for decentralization and external economies,''
\emph{Manag. Sci.}, vol. 11, no. 1, pp. 1--32, 1964.

\bibitem{KnuevenEtAl_2020}
B.~Knueven, J.~Ostrowski, and J.-P.~Watson, ``On mixed-integer programming formulations for the unit commitment problem,'' \emph{INFORMS J. Comput.}, 2020, https://doi.org//10.1287/ijoc.2019.0944.


\end{thebibliography}
\end{document}